\documentclass{amsart}
\usepackage[procnames]{listings}
\usepackage{cite}
\usepackage{color}
\usepackage{ amssymb }
\usepackage{hyperref}
\usepackage{amsmath}
\usepackage{amsmath}
\usepackage{mathtools}
\usepackage{amssymb}
\numberwithin{equation}{section}
\title{Optimal strong approximation for quadratic forms }
\author{Naser T. Sardari}
\address{Department of Mathematics, UW-Madison, Madison, WI 53706}

\email{ntalebiz@math.wisc.edu}

\date{\today}

	\newtheorem{thm}{Theorem}[section]
	
	\newtheorem{prop}[thm]{Proposition}
	
	\newtheorem{rem}[thm]{Remark}	
	\newtheorem{lem}[thm]{Lemma}
	\newtheorem{conj}[thm]{Conjecture}
	\renewcommand{\vec}[1]{\mathbf{#1}}

	\newtheorem{cor}[thm]{Corollary}
	
	\theoremstyle{defi}

	\theoremstyle{pf}

\begin{document}
\maketitle
\begin{abstract}
For a non-degenerate integral quadratic form $F(x_1, \dots , x_d)$ in $d\geq5$ variables, we prove an optimal strong approximation theorem. Let  $\Omega$ be a fixed compact subset of the affine quadric $F(x_1,\dots,x_d)=1$ over the real numbers. Take a small ball $B$  of radius $0<r<1$ inside $\Omega$, and an integer $m$.  Further assume that $N$ is a given integer which satisfies $N\gg_{\delta,\Omega}(r^{-1}m)^{4+\delta}$ for any $\delta>0$. Finally assume that an integral vector $(\lambda_1, \dots, \lambda_d) $ mod $m$ is given. Then we show that there exists an integral solution $X=(x_1,\dots,x_d)$ of $F(X)=N$ such that $x_i\equiv \lambda_i \text{ mod } m$ and $\frac{X}{\sqrt{N}}\in B$, provided that  all the local conditions are satisfied. We also show that 4 is the best possible exponent. Moreover, for a non-degenerate integral quadratic form in 4 variables we prove the same result  if $N$ is odd and $N\gg_{\delta,\Omega} (r^{-1}m)^{6+\epsilon}$. Based on our numerical experiments on the diameter of LPS Ramanujan graphs and the expected square root cancellation in a particular  sum that  appears in Remark~\ref{evidence}, we conjecture that the theorem holds for any quadratic form in 4 variables with the optimal exponent $4$.
\end{abstract}

\tableofcontents


\section{Introduction}

\subsection{Statement of results}\noindent
Before stating our main theorem, we discuss an application to a classical problem. It addresses the question of approximating  real matrices by integral matrices. More precisely, let $A=[a_{i,j}]$ be a $2\times 2$ matrix of determinant 1 and $m$ be a positive integer.  How well can one approximate $A$ by $m^{-\frac{1}{2}}H$, where $H=[h_{i,j}]$ is an integral matrix of determinant $m$? Tijdeman \cite{Tj}, showed that there exists $H$ such that 
\begin{equation*}
\max_{i,j}|m^{-\frac{1}{2}}h_{ij}-a_{ij}| <Cm^{-\frac{1}{18}}(\log m)^{7/9},
\end{equation*}
where $C$ is a  constant depending on $\max{a_{ij}}$. Later, Harman \cite{Harman} improved the exponent $-\frac{1}{18}$ to $-\frac{1}{8}$ and showed this exponent cannot be smaller than $-\frac{1}{4}$. Harman remarked that if Hooley's $R^*$ Conjecture \cite{Hoole} were true, then the exponent drops to $-\frac{1}{6}$. Subsequently, Chiu \cite[Remark 1.6]{Chiu} remarked that by assuming the $p$-th Fourier coefficient of $SL_2(\mathbb{Z})$ Maass cusp forms are less than $2p^{r}$ for all prime numbers $p$ then the exponent is less than $[(r+1/2)-1]/3$. Therefore,  if one assumes the Ramanujan Conjecture for $SL(2,\mathbb{Z})$ Maass cusp forms, then the exponent drops to $-\frac{1}{6}$.  The best known bound toward the Ramanujan conjecture is $\frac{7}{64}$ \cite{Kim} and this yields $-\frac{1}{8}-\frac{1}{192}$, which slightly improves  Harman's bound. We note that under the most favorable conditions, the two different methods give the same exponent $-\frac{1}{6}$. We show that the matrix approximation is possible with the exponent $-\frac{1}{6}$ without any assumptions. 

\begin{cor} \label{cor1}
Fix $\Omega$ a compact subset of $SL(2,\mathbb{R})$ and any $\delta>0$ . Then for every matrix $A \in \Omega$ and $m\in \mathbb{Z}$, there exists an integral matrix $H \in M_{2\times 2}[\mathbb{Z}]$ such that $\det(H)=m$ and
\begin{equation*}
\|A-\frac{1}{\sqrt{m}}H   \| \ll m^{-\frac{1}{6}+\delta},
\end{equation*}
where  $\| [a_{ij}]_{2\times 2} \|:=\sup |a_{ij}|$ and the implicit constant involved in $\ll$ only depends on $\delta$ and $\Omega$. Moreover, we cannot replace $-\frac{1}{6}$ in the exponent with a number smaller than $-\frac{1}{4}$. 
\end{cor}

We first state our main theorem in a qualitative form. Let $F(X)$ be a non-degenerate quadratic form in $d\geq 4$ variables, and  $\Omega$ be  a fixed compact subset of the affine quadric $F(X)=1$ over the real numbers. 
\begin{thm}\label{main}
Take a small ball $B$ of radius $0<r<1$ with center in $\Omega$, and an integer $m$.  Assume that  an integral vector ${\lambda}:=(\lambda_1, \dots, \lambda_d)$ mod $m$  and an integer $N$ are given such that $F(\lambda) \equiv N \mod m$. Further assume that  $F(X)=N$  has a local solution  $\vec{x}_p\in \mathbb{Z}_p^d$ for all primes $p$ such that   $\vec{x}_p\equiv \lambda \mod p^{ord_p(m)}.$ If $d\geq 5$, assume that $N$ satisfies
$$N \gg (r^{-1}m)^{4+\delta}$$
for any $\delta>0$, where the constant for $\gg$ only depends on $\Omega$ and $\delta$.  Then there exists an integral solution $F(x_1,\dots,x_d)=N $ such that
\begin{equation}
\label{sac}
x_i\equiv \lambda_i \text{ mod } m \text{ and } \frac{(x_1,\dots,x_d)}{\sqrt{N}}\in B.\end{equation}
 Moreover, the exponent $4$ in  $N \gg (r^{-1}m)^{4+\delta}$ is optimal. If $d=4$ the same result holds provided that $N$ is odd and $N \gg (r^{-1}m)^{6+\delta}$. The exponent $6$ cannot be replaced with a number smaller than $4$.
 \end{thm}
%
%
%
Based on our numerical experiments on the   diameter of LPS Ramanujan graphs \cite{Naser,Rivin, Sardaric} (see Remark~\ref{LPSrem}), and  the expected square root cancellation in a particular  sum that  appears in Remark~\ref{evidence}, we make the following conjecture.

\begin{conj}\label{mainconjecture}
We conjecture that for $d=4$ the result of Theorem~\ref{main} holds  if  $N$ is odd and 
$N\gg(r^{-1}m)^{4+\delta}.$
Therefore, the exponent $4$ is the optimal exponent for every quadratic form with $4$ or more variables.

\end{conj}
\begin{rem}
Corollary~\ref{cor1} follows from Theorem~\ref{main} when  $F(x_1, \dots,x_4)= x_1x_4-x_2x_3$. Moreover, by Conjecture~\ref{mainconjecture}, the matrix approximation is possible with the optimal exponent $-\frac{1}{4}$.
\end{rem}

For more evidence of Conjecture~\ref{mainconjecture}, we refer  the reader to the recent work of Sarnak~ \cite{Sarnak2, Sarnak3}. He proved  the existence of the optimal lift for the case of $SL_2(\mathbb{Z})\to SL_2(\mathbb{Z}/q\mathbb{Z})$ which is shown by an elementary method.

Another application of Theorem~\ref{main}  recovers the best known upper bound for  the diameter of the Ramanujan graphs first constructed explicitly by Margulis~\cite{Margulis}, and later independently by Lubotzky, Phillips and Sarnak \cite{Rama}. More precisely, let $X^{p,m}$ be the the LPS Ramanujan graph of degree $p+1$, and vertices indexed by $PGL_2(\mathbb{Z}/m\mathbb{Z})$, where $p$ is a prime number and $m\in \mathbb{Z}$; see \cite[Section 2]{Rama}.  For any $x,y \in X^{p,m}$, let $d(x,y)$ be the length of the shortest path between $x$ and $y.$ Define the diameter of $X^{p,m}$ by:
$$
\text{diam}(X^{p,m}) := \sup_{x,y \in X^{p,m}} d(x,y).
$$
It is easy to see that for  $\text{diam}(X^{p,m}) \geq \log_p|X^{p,m}|$ where $|X^{p,m}|$ is the number of vertices of $X^{p,m}.$   In \cite[Theorem 5.1]{Rama}, the authors proved that  $\text{diam}(X^{p,m}) \leq (2+\epsilon)\log_p|X^{p,m}|.$ They proved this result by appealing to the Ramanujan  bound on  the Fourier coefficients of weight two modular forms \cite{Eichler}. We recover this upper bound by applying Theorem~\ref{main} to the quadratic form $F(x_1, \dots,x_4)=x_1^2+4x_2^2+4x_3^2+4x_4^2$; see \cite[Section 3]{Rama}. In our proof, we  use Weil's bound for the Kloosterman sums.

 In fact, the diameter of the LPS Ramanujan
graphs was expected to be bounded from above by the optimal bound $(1+\epsilon)\log_p|X^{p,m}|$; see~\cite[Chapter 3]{10}. However,  Theorem~\ref{main} implies that for any prime  $p$ there exists  an infinite sequence of integers $\{m_i\}$ such that  $\text{diam}(X^{p,m_i}) \geq (4/3) \log_p|X^{p,m_i}|$; see \cite[Theorem 1.2]{Naser}.
\begin{rem}\label{LPSrem}
In  \cite{Naser,Rivin}, we demonstrated some numerical experiments on the diameter of the LPS Ramanujan graphs $X^{p,m}$, which shows that their diameter is asymptotically
$$(4/3)\log_{p}|X^{p,m}|.$$ 
In our experiments, $p$ is fixed and $m$ is growing.  This supports Conjecture~\ref{mainconjecture}. 
\end{rem}

We now state Theorem~\ref{main} in a stronger, quantitative form, which provides a lower bound on  the number of integral points inside an open neighborhood of the adelic topology. First, we define the adelic topology. 
%

For every prime  $p$, define a natural norm on $\mathbb{Q}_p^d$ via  $\|\vec{x}\|_{p}=\max_{1\leq i \leq d} |x_i|_p$ where $ |x_i|_p:=p^{-\text{ord}_p(x_i)}.$ 
For the archimedean place $\infty$, we fix an arbitrary norm $\|.\|_{\infty}$ on $\mathbb{R}^d$. For a place $v$ of $\mathbb{Q}$ ($v=p \text{ or } \infty $)  and a point $\vec{o} \in \mathbb{Q}_v^d$, we define the ball $B_{v}(\vec{o},r)$  centered at $\vec{o}$ with radius $r$ by
$
B_{v}(\vec{o},r):=\big\{\vec{x} \in \mathbb{Q}_v^d: \|\vec{x}-\vec{o}\|_v\leq r      \big\}.
$
Note that for a non-archimedean place $v$, the radius $r$ takes discrete values $q_v^n$, where $q_v$ is the order of the residue field and $n\in \mathbb{Z}$. Let $\mathbb{A}_{\mathbb{Z}}^d:= \mathbb{R}^d\times \prod_p \mathbb{Z}_p^d$ be the integral ring of adeles.   We define a global ball inside $\mathbb{A}_{\mathbb{Z}}^d$ to be a product of local ones, subject to the condition that the radius of the local balls is $1$ for all but a finite set of places.  That is,
\begin{equation}\label{BBB}
B_{\vec{a},\vec{r}}:=B_{\infty}(\vec{a}_{\infty},r)\times \prod_{p}B_p(\vec{a}_p,p^{-\nu_p}),
\end{equation}
where $\vec{a}:=(\vec{a}_{\infty}, \vec{a}_p)_p\in \mathbb{A}_{\mathbb{Z}}^d$ and  $\vec{r}:=(r,p^{-\nu_p})_p \in\mathbb{R}^{\mathbb{N}}$.  Furthermore, we are free to choose $\vec{a}_{\infty} \in \mathbb{R}^d$, $\vec{a}_p\in \mathbb{Z}_p^d$, and   $\nu_p \geq 0$ subject to the condition that $\nu_p=0$ outside a finite set of primes. By this condition,  $m:=\prod_{p}{p^{\nu_p}}$ is well-defined for some $m\in \mathbb{Z}$. We define the norm of the global ball $B_{\vec{a},\vec{r}}$ to be
$|B_{\vec{a},\vec{r}}|:=rm^{-1}.$
%
%
%
%
%
Consider the following compact topological space: 
$$\Theta:=\Omega \times \prod_{p}V_N(\mathbb{Z}_p), $$
where $V_N(R):=\{\vec{x}\in R^d: F(\vec{x})=N \}$ for any commutative ring $R.$ 
Take a global ball $B_{\vec{a},\vec{r}}$ where  $\vec{a}\in\Theta$
 and $\gcd(\vec{a}_p,p)=1$, i.e. $|\vec{a}_p|_p=1$. Define the $p$-adic density of $V_N(\mathbb{Z}_p)$ associated to the local ball $B_p(\vec{a}_p,p^{-\nu_p})$ by
$$
 \sigma_p(\vec{a}_p,p^{-\nu_p},N):=\lim_{k\to \infty }\frac{ \# \{\vec{x}\in \big(\frac{\mathbb{Z}}{p^{(k+\nu_p)}\mathbb{Z}}\big)^d: F(\vec{x})=N\mod p^{k+\nu_p} \text{ and } \vec{x}\equiv \vec{a}_p\mod p^{\nu_p}  \}}{p^{k(d-1)}}.
$$
If $\nu_p=0$ (which means we don't have any mod $p$ congruence condition) then we simply write $ \sigma_p(N)$.
It follows from  Hensel's lemma that the above limit exists.
 The condition $ \sigma_p(\vec{a}_p,p^{-\nu_p},N)\neq 0$ is equivalent to $V_N(\mathbb{Z}_p)\cap B_p(\vec{a}_p,p^{-\nu_p}) \neq \emptyset$.  Let $ \mathfrak{S}_{B_{\vec{a},\vec{r}}}(N):=\prod_{p}\sigma_{p}(\vec{a}_p,p^{-\nu_p},N)$ be the singular series. If $\nu_p=0$ for all primes $p$, then we write $ \mathfrak{S}(N):= \prod_p  \sigma_p(N) $ for the singular series. For $d\geq 4$ the product definition of $\mathfrak{S}_{B_{\vec{a},\vec{r}}}(N)$ is absolutely convergent and the singular series  is nonzero if and only if $\sigma_{p}(\vec{a}_p,p^{-\nu_p},N)\neq 0$ for every $p$; see\cite{Siegel}.  Hence, if $\mathfrak{S}_{B_{\vec{a},\vec{r}}}(N)\neq 0$ then $V_N(\mathbb{Z}_p)\cap B_p(\vec{a}_p,p^{-\nu_p}) \neq \emptyset$ for every prime $p$, and we say that all the local conditions are satisfied. 
%
%
The quantitative form of our main theorem is as follows.

\begin{thm}\label{strong}
 If $d\geq 5$, then we have 
\begin{equation}\label{formulaa}
|V_N(\mathbb{Z})\cap B_{\vec{a},\vec{r}}|\gg \mathfrak{S}_{B_{\vec{a},\vec{r}}}(N) |B_{\vec{a},\vec{r}}|^{d-1}N^{\frac{d-2}{2}}\Big(1+O\big((|B_{\vec{a},\vec{r}}|^2\sqrt{N} )^{-\frac{d-3}{2}}N^{\epsilon}\big) \Big)
\end{equation}
for some $\epsilon>0$, where the implied constant in $\gg$ and $O$ only depends on  $\epsilon$ and $\Omega$ and not on $N$ or the global ball $B_{\vec{a},\vec{r}}$. For $d=4$, we have 
\begin{equation}\label{formulaa4}
|V_N(\mathbb{Z})\cap B_{\vec{a},\vec{r}}|\gg\mathfrak{S}_{B_{\vec{a},\vec{r}}}(N) |B_{\vec{a},\vec{r}}|^{3}N\Big(1+O\big((|B_{\vec{a},\vec{r}}|^3\sqrt{N} )^{-\frac{1}{2}}N^{\epsilon}\big) \Big).
\end{equation}
\end{thm}

\begin{rem}\label{lowremark}
If  $d\geq 5$, then we  have  \cite[Remark~7 Page 73]{Malyshev} 
\begin{equation}\label{lowerdenss}c_{\epsilon}N^{-\epsilon}\leq \mathfrak{S}_{B_{\vec{a},\vec{r}}}(N),\end{equation}
where $\epsilon >0$ is any positive real number and  $c_{\epsilon}$ only depends on $\Omega$ and $\epsilon$. 
%
%
 Therefore, if $N^{\delta}\gg (|B_{\vec{a},\vec{r}}|^2\sqrt{N} )$ for some $\delta>0$ then  the error term is smaller than the main term, and as a result we have an integral point inside $B_{\vec{a},\vec{r}}$ which implies Theorem~\ref{main} for $d\geq 5.$ For $d=4$,  the inequality \eqref{lowerdenss} holds if $N$ is odd. Similarly, if  $N^{\delta}\gg (|B_{\vec{a},\vec{r}}|^3\sqrt{N} )$ then we have an integral point inside $B_{\vec{a},\vec{r}}$ and this implies Theroem~\ref{main} for $d=4$. 
%
%
%
 Morerover, if  $\gcd(p^{\nu_p},a_1(p),\dots,a_d(p))=p^{h_p}\neq1$ for some $p$, then we can divide $N$ by $p^{2h_p}$, replace the local balls at the place $p$ by $B(\frac{\vec{a}_p}{p^{h_p}},p^{-(\nu_p-h_p)})$, and use the above theorem. Since the power of $N$ is $\frac{d-2}{2}$ and the power of $|B_{\vec{a},\vec{r}}|$ is $d-1$ in the formula \eqref{formulaa},  the main term would be multiplied by $\prod_{p}p^{h_p}$. 

\end{rem}

We now state a theorem which implies that the exponent 4 in Theorem~\ref{main} is optimal.  It is based on the principle that rational points on $V_1$ of low height repel other rational points.  Let $\vec{q}:=(q_1,\dots,q_d)$ be any primitive integral vector, i.e. such that $\gcd(q_1,\dots,q_d)=1$, with $F(\vec{q})>0$. Let $\nu:=(\nu_p)$ be a sequence of nonnegative integers indexed by all primes $p$ such that $\nu_p=0$ for all but a finite set of primes, and if $\nu_p\neq 0$ then  $\Big(\frac{N}{p}  \Big)=\Big(\frac{F(\vec{q})}{p}  \Big)$  where $\Big(\frac{.}{.}\Big)$ is the  Legendre symbol. Moreover, let $a:=(a_p)$ be a sequence of $p$-adic integers $a_p\in \mathbb{Z}_p$ indexed by all prime numbers such that if $\nu_p \neq 0$, then  $a_p$ is one of the two solutions of  $a_{p}^2=\frac{N}{F(\vec{q})}$, and  if $\nu_p=0$ then $a_p=0$.  Define 
\begin{equation}\label{globalball}B_{\vec{q},\nu,a}:=B_{\infty}(a_{\infty}\vec{q},r)\times \prod_{p}B_p(a_{p}\vec{q},p^{-\nu_{p}}),\end{equation}
where  $a_{\infty}=\frac{1}{\sqrt{F(\vec{q})}}$.   Our result is the following. 
\begin{thm}\label{repthm}
 There exist positive numbers  $c_1$ and $c_2$  depending  only on $\Omega$ and $F(X)$ such that whenever one has 
\begin{equation}\label{cond11}|B_{\vec{q},\nu,a}| \geq c_1 \frac{ N^{-\frac{1}{4}}}{|\vec{q}|^{1/2}},\end{equation}
then there exists a global ball $B^{\prime}$ inside $B_{\vec{q},\nu,a}$ with no integral points and a norm greater than
$$|B^{\prime}|>c_2 \frac{ N^{-\frac{1}{4}}}{|\vec{q}|^{1/2}}.$$

\end{thm}

In the last application of our main theorem, we  give a lower bound for the number of lattice points inside a small cap on a sphere,  i.e. the intersection of a ball with a small radius and the sphere. In the appendix of \cite{Bourgain}, an upper bound for the number of lattice points in such small caps is given. Following a suggestion of Bourgain \cite[Page~23 ]{Sarnak2}, we promote this to a lower bound.\begin{cor}\label{bourgaincor}
Let  $S^{d-1}(R)$ be any $(d-1)$-dimensional sphere of radius $R$ centered at the origin  such that $N:=R^2$ is an integer and $d\geq5$. Suppose we are given a spherical cap of diameter $Y\leq R$. 
%
%
%
%
%
%
Then, the number of integral lattice points inside this spherical cap is at least 
\begin{equation*}\label{formula}
\sigma_{\infty} \mathfrak{S}(N)\frac{Y^{d-1}}{R}\Big(1+O\big( R^{\epsilon} (\frac{R}{Y^2})^{\frac{d-3}{2}}\big)\Big),
\end{equation*}
where $\sigma_{\infty}$ is a constant that only depends on $d$ and $ \mathfrak{S}(N)$ is the singular series associated to the quadric $x_1^2+\dots +x_d^2=N$. The implied constant in $O$ depends on $\epsilon$ and $d$  but not on the spherical cap or $N$. In particular, if $Y\gg_{\delta} R^{1/2+\delta}$ for any $\delta>0$ then we have an integral point inside the cap. For $d=4$, if we assume that $N$ is odd then the number of integral points is at least 
\begin{equation*}
\sigma_{\infty} \mathfrak{S}(N)\frac{Y^{3}}{R}\Big(1+O\big( R^{\epsilon} (\frac{R}{Y^{3/2}})\big)\Big).
\end{equation*}
In this case, if $Y\gg_{\delta} R^{2/3+\delta}$ then we have an integral point inside the cap.
\end{cor}

 On the other hand, it follows from Theorem~\ref{repthm} that  there are caps of diameter  $Y\gg R^{1/2}$  on any  $S^{d-1}(R)$ with no integral points inside them, where the implied constant in $\gg$ only depends on $d$. Given $R>0$ such that $R^2\in \mathbb{Z}$, we let $C(R)$ denote the maximum volume  of any cap on  $S^{d-1}(R)$ which contains  no integral points. 
%
%
Sarnak defined  \cite{Sarnak2}
 the covering exponent of integral points on the sphere by:
 \begin{equation*}
 \begin{split}
 K_d&:=\limsup_{R \to \infty}\frac{\log \big(\# S^{d-1}(R)\cap \mathbb{Z}^d \big) }{\log\big( \text{vol } S^{d-1}(R)/C(R)\big)}.
 \end{split}
 \end{equation*}
 It follows from Theorem~\ref{repthm} and Corollary~\ref{bourgaincor} that $K_d=2-\frac{2}{d-1}$ for $d\geq 5$ and $4/3\leq K_4 \leq 2$.
 
 \begin{rem}
In his letter \cite{Sarnak2} to  Aaronson and  Pollington, Sarnak showed that $4/3\leq K_4 \leq 2$. To show that $K_4\leq 2$ he appealed to the Ramanujan bound on the Fourier coefficients of weight $k$ modular forms, while the lower bound $4/3\leq K_4$ is a consequence of an elementary number theory argument. Furthermore, Sarnak states some open problems \cite[Page $24$]{Sarnak2}.  The first one is to show that  $K_4<2$ or  even that $K_4=4/3$. Conjecture~\ref{mainconjecture} implies  that 
$K_4=4/3$.
 \end{rem}

More generally, Ghosh, Gorodnick and Nevo studied the covering exponent of the orbits of a lattice subgroup $\Gamma$ in a connected Lie group $G$, acting on a suitable homogeneous spaces $G/H$; see \cite{GGN, GGN1,GGN2}. They linked  the covering exponent of $\Gamma$ to the spectrum of $H$ in the automorphic representation on $L^2(\Gamma\backslash G).$ In particular, they showed that $K_{\Gamma} \leq 2$ if  the restriction of the unitary representation on $L^2(\Gamma\backslash G)$ to $H$ has tempered spherical spectrum as a representation of $H$; see\cite[Theorem 3.5]{GGN2} and \cite[Theorem 3.3]{GGN1}.  This recovers  the above result of Sarnak for $d=4.$  For $d\geq 5$, by using the best bound on  the generalized Ramanujan conjecture for $SO_{d},$ they showed that   \cite[Page 12]{GGN}  $1 \leq K_d\leq 4-4/(d-1) $ for odd $d$ and $1 \leq K_d\leq  4-16/(d+2)$ for even $d$. They raised the question of improving  these bounds in \cite[Page 11]{GGN}. As pointed out above, we give a definite answer to this question and show that $K_d=2-\frac{2}{d-1}$ for $d\geq 5.$

\subsection{Further motivations and techniques }\noindent

In several papers, Wright proved various results about the representation of an integer  $N$ as a sum of squares of integers ``almost proportional" to assigned positive real  numbers $\lambda_1, \lambda_2, \dots, \lambda_d$. 
 In \cite{Wright}, he showed that if   $N^{1-\frac{1}{8}+\epsilon} \leq U$ for some $0 < \epsilon$ and $
\lambda_1+\dots+\lambda_5=1,  \text{ where }  0 < \lambda_i,$
then there exists an integral solution $(n_1, \dots, n_5)$ to $
N=n_1^2+\dots+n_5^2,$ where
$|n_i^2 -\lambda_i N|< U.$
Note that by Theorem~\ref{main}  the inequality $N^{1-\frac{1}{8}+\epsilon} \leq U$ is not sharp and can be improved to $N^{1-\frac{1}{4}+\epsilon} \leq U$. He also showed that the number of representation is $
\asymp \frac{U^4}{\sqrt{N}}.
$
By an entirely elementary method Auluck and Chowla   in \cite{Chowla} proved a sharp result for the special point $(\lambda_1,\dots,\lambda_4)=(\frac{1}{4},\frac{1}{4},\frac{1}{4},\frac{1}{4})$ and sum of four squares. They showed that  every positive integer $N\neq 0$ mod 8 is expressible in the form $N=n_1^2 +\dots+n_4^2,$
where $n_i$ are integers satisfying
$\frac{N}{4}-n_i^2=O(N^{\frac{3}{4}}).$
It follows from further work of Wright \cite{Wright2} that  $O(N^{\frac{3}{4}})$ in this theorem cannot be replaced by $o(N^{\frac{3}{4}})$. More generally, he considered the sum of $d$ $k$th powers and proved that there exists an infinite sequence of integers $\{N_i\} $ such that the diagonal point $W_{(d,k,N_i)}:=(a_i,\dots,a_i)$, where $da_i^k=N_i$ for some $a_i\in\mathbb{R}$, repels the integral points on $x_1^k+\dots+x_d^k=N_i.$ More precisely, he showed that the ball $B_{W_{(d,k,N_i)},\gamma_{(d,k)}N_i^{1/2k}}$ centered at $W_{(d,k,N_i)}$  with radius $\gamma_{(d,k)}N_i^{1/2k}$ for some fixed $\gamma_{(d,k)}>0$ contains no integral points $(x_1,\dots,x_d)$ such that  $x_1^k+\dots+x_d^k=N_i.$ We will discuss this repulsion property for $k = 2$ in more detail in Section~\ref{repel}.  In a recent paper \cite{Lower}, Daemen gave a lower bound for the number of  integral points $(x_1,\dots,x_d)$ close to the diagonal point $W_{(d,k,N)}$ such that
$x_1^k+\dots+x_d^k=N.$
  He proved that for every $k$ there exists $d_k$  such that if  $d \geq d_k$, then one has a lower bound for the number of integral points inside $B_{W_{(d,k,N_i)},N^{1/2k+\epsilon}}$ for any $\epsilon>0$ and large enough  $N$.  This lower bound differs by a bounded scalar compared to the upper. For the sum of squares, his result implies that   $d_2 \leq 9$.
He remarked that by working a little harder, one can prove that $d_2\leq 7$. Theorem~\ref{main} implies  that $ d_2\leq 5$.





We shall prove Theorem~\ref{main} using a version of the delta method, which in turn is based on the Kloosterman circle method. This method was  developed by Kloosterman   to prove the local-to-global principle  for quadratic forms in $4$ variables;  see~\cite{Kloos}. For the purpose of 
representing integers by quadratic forms,
if the number of variables is 5 or more the classical 
circle method of ``major and minor" arcs works fine.
However, for 4 variables it does not work. 
Kloosterman introduced a new method of dissection by the Farey sequence 
(with no minor arcs) which deals with 4 variables.

In our work, the method proves
to be decisive in 5 or more variables as it 
gives the optimal exponent for lifting solutions. Our optimal result for $5$ variables depends crucially on obtaining square root cancellation in certain exponential sums, and also relies on a refinement of the Kloosterman method developed by Duke, Friedlander and Iwaniec known as the delta method; see~\cite{delta}. In the delta method we use a smooth cut off function over the Farey dissections. The delta method allows us to use the rapid decay of the Fourier transform of the weight function and makes computations handier.  For proving  Theorem~\ref{main}, we use a version of this method developed by Heath-Brown in \cite{Heat}. In that method we first apply the delta method and then a Poisson summation over the sum of lattice points. Our main innovation is to introduce a special coordinate system in Section~\ref{weightfunction} which is crucial in improving the current bounds on the oscillatory integrals appearing in the delta method.

The problem of the distribution of integral points on quadrics between different residue classes to a fixed integer $m$ has been studied by Malyshev; see  \cite{Malyshev}.
Malyshev used Kloosterman method and proved a result about the distribution of integral points of the quadric $F(X)=N$ with 4 or more variables between residue class of an integer $m$. An application of Theorem~\ref{main} significantly improves the exponents of $m$ and $N$ in his main theorem for $F(X)$ fixed.  

\subsection*{Acknowledgements}\noindent
I would like to thank my Ph.D. advisor Peter Sarnak for several insightful and inspiring conversations during the course of this work. In fact the starting point of this work was his letter \cite{Sarnak2} and in particular the key observation of Bourgain noted there as well as above.  I would like to thank Professor Jianya Liu and Professor Tim Browning for their comments on the earlier versions of this work. Finally, I am grateful  for  the  comments of Simon Marshall, Masoud Zargar, and the anonymous referees that   improved the writing of this manuscript.

\section{Repulsion of integral points }\label{repel}
In this section we prove Theorem~\ref{repthm}. We assume that $F(X)=X^TAX$ is a non-degenerate quadratic form. Recall the definition of $B_{\vec{q},\nu,a}$ in \eqref{globalball}, and the notation used while formulating Theorem~\ref{repthm}. 
%
%
%

%
\begin{proof}Let $B_{\vec{q},\nu,a}^{\infty}:=\prod_{p}B_p(a_{p}\vec{q},p^{-\nu_p})$ be the finite part of the global ball $B_{\vec{q},\nu,a},$ and $m:=\prod_p p^{\nu_p}.$
We assume that $X\in V_N(\mathbb{Z})$ is an integral point such that $X\in B_{\vec{q},\nu,a}^{\infty}$. Hence, $X\equiv a_{p}\vec{q}   \text{ mod } p^{\text{ord}_p(m)}.$  
By the Chinese remainder theorem there exists $\alpha\in \mathbb{Z}$ such that
$X\equiv \alpha \vec{q}  \text{ mod } m,$
where $\alpha \equiv a_{p}  \text{ mod } p^{\text{ord}_p(m)} $. Hence, $X=m\vec{t}+\alpha \vec{q}$ for some integral vector $\vec{t}$.  We write the Taylor expansion of $F(mt+\alpha q)$ at $\alpha q$  , i.e.
$F(m\vec{t}+\alpha \vec{q})=m^2F(\vec{t})+2m\vec{t}^{T}A\alpha \vec{q}+F(\alpha \vec{q})=N.$
Since $\gcd(\alpha,m)=1$, we deduce that
$\vec{t}^{T}A \vec{q} \equiv \frac{N-F(\alpha \vec{q})}{2m\alpha} \text{  mod }  m.$
Therefore, there exists a fixed number $l_0$ mod $m^2$ such that 
\begin{equation}\label{intcon}
\left<X,A\vec{q}   \right> \equiv \left<m\vec{t}+\alpha \vec{q},A\vec{q}   \right> \equiv l_0  \text { mod } m^2.
\end{equation}
Now consider the infinite part of the global ball $B_{\vec{q},\nu,a}$, namely 
$B_{\infty}:=B_{\infty}(a_{\infty}\vec{q},r).$ Without loss of generality we  assume that $B_{\infty}$ intersects only one connected component of $V_1(\mathbb{R}).$
 We wish to find a real point $\tilde{\vec{q}} \in \sqrt{N}B_{\infty} \cap V_N(\mathbb{R})$ such that
\begin{equation}\label{conn}
\left<\tilde{\vec{q}}, A\vec{q} \right>=\frac{m^2}{2}+l_0+km^2   \text  {  for some  } k\in \mathbb{Z},  \text{ and } |\tilde{\vec{q}}-\sqrt{N}a_{\infty}\vec{q}| \text{ is minimal} .
\end{equation}
We will deduce the existence of $\tilde{\vec{q}}$ from  assumption~\eqref{cond11} (although the $\tilde{\vec{q}}$ we produce is not necessarily unique).  By the connectivity  of $\sqrt{N}B_{\infty} \cap V_N(\mathbb{R})$ and the continuity of the inner product, it suffices to show that the length of the interval $I:=\{\left<\vec{h}, A\vec{q} \right> \in \mathbb{R}: \vec{h} \in  \sqrt{N}B_{\infty} \cap V_N(\mathbb{R})  \}$ is bigger than $m^2.$  For any $r<1$ and $\vec{y}\in \Omega$, let $J_{\vec{y},r}:=\{\left<\vec{l},A\vec{y} \right>: \vec{l} \in V_1(\mathbb{R}), |\vec{l}-\vec{y}|< r \}\subset \mathbb{R}.$ Since $\Omega$ is compact and $F(X)$ is non-degenerate,  it follows that there exists a constant $c_3$ depending on $\Omega $ such that the length of $J_{\vec{y},r}$ is bigger than $c_3r^2$, independently of $\vec{y}\in \Omega.$ Moreover, there exist constants $c_4$ and $c_5$ depending on $\Omega$ such that $ c_4|\vec{q}|\leq  \sqrt{F(\vec{q})}\leq c_5|\vec{q}|$. This implies that $|J_{a_{\infty}\vec{q},r}| \geq c_3 r^2$, or equivalently that $|I|\geq c_3 r^2 \sqrt{N |F(\vec{q})|} \geq c_3 c_4 r^2 \sqrt{N} |\vec{q}|.$ By assumption~\eqref{cond11}, we have  $m^2 \leq c_1  r^{2}|\vec{q}|\sqrt{N}$. By choosing $c_1<c_3c_4$, it follows that $|I|\geq m^2$, and hence there exists $\tilde{\vec{q}}$ satisfying condition~\eqref{conn}.

Next, we give an upper bound on $|\tilde{\vec{q}}-\sqrt{N}a_{\infty}\vec{q}|$. By condition~\eqref{conn} and the connectivity of  $ \sqrt{N}B_{\infty} \cap V_N(\mathbb{R})$, it follows  that
$
|\left<\tilde{\vec{q}}-\sqrt{N}a_{\infty}\vec{q},A\vec{q} \right>|\leq m^2.
$ We write $\tilde{\vec{q}}-\sqrt{N}a_{\infty}\vec{q}=\vec{q}_1+ \vec{q}_2$, where $\vec{q}_1$ is parallel to $A\vec{q}$ and $\vec{q}_2$ is orthogonal to $A\vec{q}.$ We have
\[
|\vec{q}_1|= \frac{|\left<\vec{q}_1,A\vec{q} \right>|}{|A\vec{q}|}= \frac{|\left<\tilde{\vec{q}}-\sqrt{N}a_{\infty}\vec{q},A\vec{q} \right>|}{|A\vec{q}|}\leq \frac{m^2}{|A\vec{q}|}\ll_{\Omega}  \frac{m^2}{|\vec{q}|}.
\]
Since $F(\tilde{\vec{q}})= F(\sqrt{N}a_{\infty}\vec{q} )=N$, we have
\[
|\vec{q}_2^T A\vec{q}_2|= |\vec{q}_1^T A\vec{q}_1+2\vec{q_1}A\sqrt{N}a_{\infty}\vec{q}| \ll |\vec{q}_1|^2+ m^2\sqrt{N}a_{\infty}. 
\]
Hence, $|\vec{q}_2| \ll |\vec{q}_1|+m N^{1/4} a_{\infty}^{1/2}. $ By the triangle inequality and the assumption~\eqref{cond11},
\begin{equation}\label{residualb}
|\tilde{\vec{q}}-\sqrt{N}a_{\infty}\vec{q}|=|\vec{q}_1+ \vec{q}_2| \ll \frac{m^2}{|\vec{q}|}+m N^{1/4} a_{\infty}^{1/2} \ll m N^{1/4} a_{\infty}^{1/2}.
\end{equation}
%

Assume that $X=\tilde{\vec{q}}+\vec{\xi}$ is an integral point inside  $B_{\vec{q},\nu,a}^{\infty}$ with $F(X)=N$. We write the Taylor expansion of $F(X)$ at $\tilde{\vec{q}}$ and obtain
$ F(X)= F(\tilde{\vec{q}})+2 \tilde{\vec{q}}^TA\vec{\xi} +F(\vec{\xi}).  $
Since $F(X)=F(\tilde{\vec{q}})=N$, we deduce that 
$2 \tilde{\vec{q}}^TA\vec{\xi} +F(\vec{\xi})=0.$
Hence,
 \begin{equation*}
 |F(\xi)| +|\left< 2\xi,A (\tilde{\vec{q}}-\sqrt{N}a_{\infty}\vec{q} )\right>| \geq |\left< 2\vec{\xi},\sqrt{N}a_{\infty}A \vec{q} \right>|.
 \end{equation*}
By~\eqref{intcon} and the definition of $\tilde{\vec{q}}$  in (\ref{conn}), $ |\left< 2\vec{\xi},\sqrt{N}a_{\infty}A \vec{q} \right>| \geq m^2a_{\infty}\sqrt{N}.$ By the Cauchy-Schwarz inequality and inequality~\eqref{residualb},
$$|\left< 2\xi,A (\tilde{\vec{q}}-\sqrt{N}a_{\infty}\vec{q} )\right>| \leq 2|\vec{\xi}||A ||(\tilde{\vec{q}}-\sqrt{N}a_{\infty}\vec{q} )| \ll m N^{1/4} a_{\infty}^{1/2}|\vec{\xi}|.$$
 Therefore,
\begin{equation*}\label{dada}
|F(\vec{\xi})| +m N^{1/4} a_{\infty}^{1/2}|\vec{\xi}| \gg  m^2a_{\infty}\sqrt{N}. 
\end{equation*}
%
%
We deduce that $|\vec{\xi}| \gg   ma_{\infty}^{1/2}N^{1/4}.$
Since $a_{\infty}=1/\sqrt{F(\vec{q})}\gg \frac{1}{|\vec{q}|}$, then  there exists a constant $c_2$ depending on $\Omega$ and $F$ such that 
$\frac{|\xi|}{\sqrt{N}} \geq c_2 \frac{mN^{-1/4}}{|\vec{q}|^{1/2}}.$ Let $r^{\prime}:=c_2 \frac{mN^{-1/4}}{|\vec{q}|^{1/2}}$ and $
B^{\prime}:=B(\frac{\tilde{\vec{q}}}{\sqrt{N}},r^{\prime})\prod_{p}B(a_{p}\vec{q},p^{-\text{ord}_p(m)}). 
$
Then, $
|B^{\prime}|=c_2 \frac{mN^{-1/4}}{|\vec{q}|^{1/2}}
$ and $B^{\prime}$  does not contain any integral point. This concludes the proof of Theorem~\ref{repthm}.
%
\end{proof}

\section{The delta method}\label{sdelta}

We now begin the proof of Theorem~\ref{strong}.  In this section, we define a smooth sum $N(w, \lambda)$ that gives a lower bound for the number of integral points we wish to study, and use the delta method to give an expansion for $N(w,\lambda)$.

Assume that $F(X)=X^TAX$ is a non-degenerate quadratic form. Recall the notations used while formulating  Theorem~\ref{strong} for $F(X)$. Assume that $X:=(x_1,\dots,x_d)$ is an integral point inside the given global ball $B_{\vec{a},\vec{r}}$. Hence,
\begin{equation}\label{syscong}x_i\equiv a_i(p)   \text{ mod } p^{\nu_p},   \text{  for every } p,\end{equation}
where $\vec{a}_p=\big(a_1(p),\dots,a_d(p)\big) \in \mathbb{Z}_p^d.$ Recall that $|B_{\vec{a},\vec{r}}|=m^{-1}r$ where $m:=\prod_{p}p^{\nu_p}$. By the Chinese remainder theorem there exists $\lambda=(\lambda_1,\dots,\lambda_d) \mod m$ such that
\begin{equation}\label{modm}x_i\equiv \lambda_i  \text{ mod } m,\end{equation}
if and only if the system of congruence conditions \eqref{syscong} hold. So, counting integral points inside $B_{\vec{a},\vec{r}}$ is the same as counting integral points inside the ball $B_{\infty}(\vec{a}_{\infty},r)\subset\mathbb{R}^d$ subjected to the congruence condition~\eqref{modm}.

Next, we give an expansion of the delta function which is developed by Duke, Friedlander and Iwaniec \cite{delta}. 
 We cite this theorem from \cite[Theorem~1]{Heat}.
\begin{thm}\label{Heat}
For any integer n let
\begin{equation*}
\delta(n) =\begin{cases} 1   \text{  if  }  n=0, \\  0  \text{   otherwise.  }     \end{cases}
\end{equation*}
Then for any $Q>1$ there is a positive constant $C_{Q}$, and a smooth function $h(x,y)$ defined on the set $(0,\infty)\times \mathbb{R}$, such that
\begin{equation}\label{delta}
\delta(n)=\frac{C_{Q}}{Q^{2}}\sum_{q=1}^{\infty} {\sum_{a}}^*e(\frac{an}{q})h(\frac{q}{Q},\frac{n}{Q^2}),
\end{equation}
where ${\sum}^*$ means the sum is over $a \mod q$ with $\gcd(a,q)=1$, and the constant $C_Q$ satisfies
$
C_{Q}=1+O_N(Q^{-N}),
$
for any $N>0$. Moreover $h(x,y)\ll x^{-1}$ for all $y$, and $h(x,y)$ is non-zero only for $x \leq \max(1,2|y|)$.
\end{thm}

Let $w$ be a smooth compactly supported weight function defined on $\mathbb{R}^n$ such that 
\begin{equation}\label{property1}  w(\vec{x})=0  \text{ if }   \vec{x}\notin B_{\infty}(\vec{a}_{\infty},r).\end{equation}
%
%
%
Assume that $X\in \mathbb{Z}^d$ satisfies the condition~\eqref{modm}. We uniquely write $X=m\vec{t}+ \lambda,$ where $\vec{t}\in  \mathbb{Z}^d$ and $\lambda=(\lambda_1,\dots,\lambda_d)$ for $-\frac{m}{2} < \lambda_i \leq \frac{m}{2}$. Define 
\begin{equation}\label{defk}
k:= \frac{N-F(\lambda)}{m}.
\end{equation}
Since
$F(X) =N,$ then
$m^2F(\vec{t})+2m\lambda^T A\vec{t} =N-F(\lambda)$
which implies   $m|2\lambda^TA\vec{t}-k.$ 
 Then, 
$F(\vec{t})+\frac{1}{m}(2\lambda^T A\vec{t}-k )=0.$
We also define $$G(\vec{t}):=\frac{F(m\vec{t}+\lambda)-N}{m^2}=F(\vec{t})+\frac{1}{m}(2\lambda^T A\vec{t}-k ).$$
Finally, we  define $$N(w,\lambda):=\sum_{\vec{t}} w(m\vec{t}+\lambda) \delta{(G(\vec{t}))},$$  where $\vec{t}\in \mathbb{Z}^{d}$. Note that $N(w,\lambda)$ is the weighted number of integral points $X$ satisfying condition~\eqref{modm}. We wish to apply the delta expansion in \eqref{delta} to   $ \delta{(G(\vec{t}))}.$ Note that \eqref{delta} holds only for $n\in\mathbb{Z}.$ Moreover,  $G(\vec{t})\in \mathbb{Z}$ if and only if $m|2\lambda^TA\vec{t}-k.$ Hence, we write
\begin{equation*}
N(w,\lambda)=\frac{1}{m}\sum_{l}\sum_{\vec{t}} e\big(\frac{l}{m} (2\lambda^T A\vec{t}-k)  \big) w(m\vec{t}+\lambda) \delta{(G(\vec{t}))},
\end{equation*}
where $l$ varies mod $m.$ Then, we apply \eqref{delta} with 
 $Q:=\frac{\sqrt{N}}{mr^{-1}}$ and obtain

\begin{equation*}
\begin{split}
N(w,\lambda)=
\frac{C_Q  } {mQ^2} \sum_{l} \sum_{q} {\sum_{a}}^* \sum_{\vec{t} } 
e\Big(\frac{(lq+a)(2\lambda^T A\vec{t} - k) + am F(\vec{t})}  {mq}  \Big)
\\ 
h\Big(\frac{q}{Q},\frac{G(\vec{t})}{Q^2}\Big) w(m\vec{t}+\lambda) .
\end{split}
\end{equation*}
We apply the Poisson summation formula on the $\vec{t}$ variable to obtain
\begin{equation}\label{form13}
N(w,\lambda)=\frac{C_Q  } {mQ^2} \sum_{l} \sum_{q} {\sum_{a}}^* \sum_{\vec{c} }(mq)^{-d} S_{m,q}(a,l,\vec{c})I_{m,q}(\vec{c}),
\end{equation}
where $\vec{c}\in \mathbb{Z}^d$, $l \in (\mathbb{Z}/m\mathbb{Z})$, and $a \in (\mathbb{Z}/q\mathbb{Z})^*$, and $I_{m,q}$ and $S_{m,q}$ are given by 
\begin{equation}\label{intt}
I_{m,q}(\vec{c}):= \int h(\frac{q}{Q},\frac{G(\vec{t})}{Q^2}) w(m\vec{t}+\lambda) e\big(-\frac{\left< \vec{c},\vec{t} \right>}{mq}    \big) dt_1\dots dt_d,
\end{equation}
\begin{equation}\label{salam}
S_{m,q}(a,l,\vec{c}):=\sum_{\vec{b}\in (\mathbb{Z}/mq \mathbb{Z} )^d} e \Big(\frac{(lq+a)(2\lambda^T A\vec{b}-k) + am F(\vec{b}) + \left< \vec{c},\vec{b} \right>} {mq}\Big).
\end{equation}
We define the sum $S_{m,q}(\vec{c})$ as:
\begin{equation}\label{sala}S_{m,q}(\vec{c}):=\sum_{l }{\sum_{a}}^* S_{m,q}(a,l,\vec{c}).\end{equation}
 Therefore
\begin{equation}\label{newequu}N(w,\lambda)=\frac{C_Q}{mQ^2}\sum_{q}\sum_{\vec{c}}(mq)^{-d}S_{m,q}(\vec{c})I_{m,q}(\vec{c}).\end{equation}

\noindent Since $Q=\frac{\sqrt{N}}{mr^{-1}}$, in the Section~\ref{weightfunction}, we construct a smooth function $w$ with compact support such that
\begin{equation}\label{property2}
 w(m\vec{t}+\lambda)\neq 0 \text{ only if } \frac{G(\vec{t})}{Q^2} \ll 1.
\end{equation}
Since the support of the smooth function $h(x,y)$ is inside the interval $0 \leq x\leq \max(2y,1)$. Hence,  the summation on $q$ in formula~(\ref{newequu}) is restricted to $1\leq q \ll Q$.

%
%
%
%
%


\subsection{ The properties of the smooth function $h(x,y)$}\label{h}
In this section we cite the basic properties of the smooth function $h(x,y)$ from \cite{Heat}. 
\begin{lem}\label{L12bound}
We have
$$\int \Big|x^ky^l \frac{\partial^k  h(x,y)}{\partial x^k }   \Big| dy \ll_{l} x^l,$$
for any $l\geq 0$ and $k\in\{0,1\}.$
\end{lem}
\begin{proof} This is an easy consequence of \cite[Lemma 5]{Heat}.
\\
\end{proof}

The following lemma shows $h(x,y)$ converges to the delta function rapidly as $x\to0$. For the proof we refer the reader to \cite[Lemma 9]{Heat}.
\begin{lem}\label{del}
Let $f$ be a smooth function with compact support. Then if $x\ll 1$ we have
\begin{equation}
\int f(y)h(x,y) dy= f(0)+O_{M,f}(x^M).
\end{equation}
\end{lem}

The following lemma shows the smooth property of the $h(x,y)$ in the $y$ variable. Since it converges to the delta function from the previous lemma, as $x\to 0$ the decay rate of the Fourier transform is slower as $x\to 0$. For the proof, we refer the reader to \cite[Lemma 17]{Heat}.
\begin{lem}\label{Fourier}
Let $w$ be a smooth compactly supported weight function and $p(t,x)$ be the Fourier transform of $w(y)h(x,y)$ in the $y$ variable, i.e., 
\begin{equation*}
p(t,x)=\int_{-\infty}^{\infty}  w(y)h(x,y)e(-ty) dy.
\end{equation*}
Then $p(t,x)$ decays faster than any polynomial in the variable $xt$, i.e. for every $N\geq 0$ we have
\begin{equation*}
p(t,x) \ll _{w, N } (xt)^{-N}.
\end{equation*}
\end{lem}


\section{Quadratic exponential sums $S_{m,q}(\vec{c})$}\label{mm}

Recall the definition of $S_{m,q}(\vec{c})$ from \eqref{sala}. In this section, we prove the following upper bound on the average norm of $S_{m,q}(\vec{c})$ with respect to $q$. A version  of this inequality was proved  in Heath-Brown's paper \cite[Lemma 28]{Heat}. Our proof uses Weil's bound (resp. Sali\'e's bound) on generalized Kloostermann sums (resp. Sali\'e sum) for even (resp. odd) dimensions.

\begin{prop}\label{mM}
We have the following upper bound 

\begin{equation*}\label{weil}
\sum_{q=1}^{X}m^{-d}q^{-\frac{d+1}{2}}|S_{m,q}(\vec{c})|   \ll_{\Delta} m^{\epsilon}X^{1+\epsilon},
\end{equation*}
where $X=O(N^{A})$ for any fixed power $A$.

\end{prop}

We give the proof of Proposition~\ref{mM} at the end of this section. We begin by proving some auxiliary lemmas. 

\begin{lem}\label{conglem} Unless $c\equiv- 2(lq+a)A\lambda  \text{ mod } m$,
we have
 $S_{m,q}(a,l,\vec{c})=0$. As a result $S_{m,q}(\vec{c})=0$, unless 
$\vec{c}\equiv \alpha A\lambda  \text{ for  some scalar $\alpha$ mod $m$}$.

\end{lem}
\begin{proof}We write the  vector $\vec{b}$ in the sums $S_{m,q}(a,l,\vec{c})$ as 

$$\vec{b}=q\vec{b}_1+\vec{b}_2,$$
where $\vec{b}_1$ is a vector mod $m$ and $\vec{b}_2$ is a vector mod $q$. Then we write $S_{m,q}(a,l,\vec{c})$  as a summation over $\vec{b}_1$ and $\vec{b}_2$, and  obtain
\begin{equation*}\label{ssss}
S_{m,q}(a,l,\vec{c})=\sum_{\vec{b}_2} e\Big(\frac{(lq+a)(2\lambda^TA\vec{b}_2-k)+amF(\vec{b}_2)+\left<\vec{c},\vec{b}_2  \right>}{mq}   \Big)  \sum_{\vec{b}_1} e\Big(  \frac{(lq+a)2\lambda^TA\vec{b}_1+\left< \vec{c},\vec{b}_1 \right>}{m}   \Big).
\end{equation*}
It is easy to check that the summation over $\vec{b}_1$ is zero, unless $\vec{c}\equiv-2(lq+a)A\lambda \text{ mod } m.$  Hence, we conclude the lemma.

 \end{proof}

\noindent In the rest of this section we give an upper bound on $S_{m,q}(\vec{c})$. Let $\Delta:=\det A$. By~\eqref{sala}
$$S_{m,q}(\vec{c})=\sum_{l }{\sum_{a }}^* \sum_{\vec{b}} e \Big(\frac{(lq+a)(2\lambda^TA\vec{b}-k) + amF(\vec{b}) + \left< \vec{c},\vec{b} \right>} {mq}\Big).$$
Since the summation over $l$ is  nonzero only if $m|2\lambda^TA\vec{b}-k$ and it is $m$ when $m|2\lambda^TA\vec{b}-k$,  we have
\begin{equation}\label{ssume}
S_{m,q}(\vec{c})=m\sum_{\vec{b}, a} e \Big(\frac{a(2\lambda^TA\vec{b}-k) + amF(\vec{b}) + \left< \vec{c},\vec{b} \right>} {mq}\Big),
\end{equation} 
where the summation is over  $a\in\big(\mathbb{Z}/q\mathbb{Z}\big)^*$ and vectors $\vec{b}\in\big(\mathbb{Z}/mq\mathbb{Z}\big)^d$ such that $m|2\lambda^TA\vec{b}-k.$
We assume that $q=q_1q_2,$
where $\gcd(q_1,2\Delta m)=1$ and the set of primes which divide $q_2$ is a subset of prime divisors  of $2\Delta m$. So, $\gcd(q_1,mq_2)=1.$ By the Chinese remainder theorem  we  write \begin{equation}\label{decomk}k=mq_2k_1+q_1k_2\end{equation} 
for some integers $k_1$ and $k_2$, and 
$$a=q_2a_1+q_1a_2,$$
where $a_1\in \big(\mathbb{Z}/q_1\mathbb{Z}\big)^*$ and  $a_2\in \big(\mathbb{Z}/q_2\mathbb{Z}\big)^*.$  We also write 
$$
\vec{b}=mq_2 \vec{b}_1 +q_1\vec{b}_2,
$$
where   $\vec{b}_1 \in \big(\mathbb{Z}/q_1\mathbb{Z}\big)^d$ and  $\vec{b}_2 \in \big(\mathbb{Z}/mq_2\mathbb{Z}\big)^d$ such that $m|2\lambda^TA\vec{b}-k$.  
We substitute  $k=mq_2k_1+q_1k_2,$ $a=q_2a_1+q_1a_2$, and $\vec{b}=mq_2 \vec{b}_1 +q_1\vec{b}_2$ in formula (\ref{ssume}) for $S_{m,q}(\vec{c})$, and obtain
\begin{equation}\label{FF}
S_{m,q}(\vec{c})=S_1S_2,
\end{equation}
where
\begin{equation}\label{S1}
S_1:=\sum_{a_1,\vec{b}_1} e\Big(\frac{2q_2a_1\lambda^TA\vec{b}_1+a_1(mq_2)^2F(\vec{b}_1)+\left<\vec{c},\vec{b}_1\right>-q_2a_1k_1    }{q_1}               \Big),
\end{equation}
and 
\begin{equation}\label{S2}
S_2:=m\sum_{a_2,\vec{b}_2}  e\Big(\frac{2q_1a_2\lambda^TA\vec{b}_2+a_2mq_1^2F(\vec{b}_2)+\left<\vec{c},\vec{b}_2\right>-q_1a_2k_2}{mq_2}     \Big).
\end{equation}
 In the following lemma we give upper bounds for $S_1$. This lemma and its argument  is similar to \cite[Lemma 26]{Heat}.

\begin{lem}\label{1}
Let $G(X):=X^TBX$ and $D:=\det(B)\neq 0$, where $B$ is a symmetric $d\times d$ integral matrix. Let $q$ be an integer such that $\gcd(q,2D)=1$.  For $t\in \mathbb{Z}/q\mathbb{Z}$ and $\vec{c},\vec{c}^{\prime}\in \big(\mathbb{Z}/q\mathbb{Z}\big)^d $, we define
\begin{equation*}
S(G,\vec{c},\vec{c}^{\prime},t):=\sum_{a,\vec{b}} e\Big(\frac{a(G(\vec{b})+\left<\vec{c}^{\prime},\vec{b}\right>+t)+\left<\vec{c},\vec{b}\right>    }{q}               \Big),
\end{equation*}
where the sum being taken over all  $a\in\big(\mathbb{Z}/q\mathbb{Z}\big)^{*}$ and $\vec{b}\in\big(\mathbb{Z}/q\mathbb{Z}\big)^d$. Then
\begin{equation}\label{formulaaw}
S(G,\vec{c},\vec{c}^{\prime},t)=\Big(\frac{D}{q}\Big)\tau_{q}^d Kl(G,\vec{c},\vec{c}^{\prime},t),
\end{equation}
where $\tau_{q}:=\sum_{x}e\big(\frac{x^2}{q} \big)$ is the Gauss sum, $\Big(\frac{.}{.}\Big)$ is the Jacobi symbol, and $Kl(G,\vec{c},\vec{c}^{\prime},t)$ is either a Kloosterman sum or a Sali\'e sum. As a result, we have
\begin{equation*}
|S_1|\leq q_{1}^{\frac{d+1}{2}}\tau(q_1)\gcd(q_1,N)^{1/2}.
\end{equation*}
\end{lem}
\begin{rem}
We note that the exponent $\frac{d+1}{2}$ in the above upper bound is optimal, and corresponds to the full square root cancellation in both the $a$ and $\vec{b}$ variables.
\end{rem}
\begin{proof}Since $q$ is odd we can diagonalize our quadratic form $G(X)$ mod $q$, and write $G(X)=\sum_{i=1}^{d} \alpha_i x_i^2.$
Hence,
$$ S(G,\vec{c},\vec{c}^{\prime},t):={\sum_{a}}^* e\Big(\frac{at}{q}\Big)  \prod_{j=1}^{d} \sum_{b} e\Big(\frac{a(\alpha_jb^2+c_j^{\prime}b)+c_jb}{q}   \Big),  $$
where $b\in \mathbb{Z}/q\mathbb{Z}.$
We complete the square to obtain
\begin{equation*}\begin{split} S(G,\vec{c},\vec{c}^{\prime},t):&={\sum_{a}}^* e\Big(\frac{at}{q}\Big)\prod_{j=1}^{d} \sum_{b} e\Big(\frac{a\alpha_j\big(b+\frac{ac_j^{\prime}+c_j}{2a\alpha_j}\big)^2-\frac{(ac_j^{\prime}+c_j)^2}{4a\alpha_j}}{q}   \Big)
\\
 &={\sum_{a}}^* e\Big(\frac{at}{q}\Big)\prod_{j=1}^{d} e\Big( \frac{-\frac{(ac_j^{\prime}+c_j)^2}{4a\alpha_j} }{q} \Big) \sum_{b} e\Big(\frac{a\alpha_j\big(b+\frac{ac_j^{\prime}+c_j}{2a\alpha_j}\big)^2}{q}   \Big).  \end{split}\end{equation*}
We note that 
$$\sum_{b} e\Big(\frac{a\alpha_j\big(b+\frac{ac_j^{\prime}+c_j}{2a\alpha_j}\big)^2}{q}   \Big)= \Big(\frac{a\alpha_j}{q}\Big)\tau_q,$$
where $\tau_{q}:=\sum_{x}e\big(\frac{x^2}{q} \big)$ is a quadratic Gauss sum and  $\Big(\frac{a\alpha_j}{q}\Big)$ is the Jacobi symbol. Hence,
$$S(G,\vec{c},\vec{c}^{\prime},t)=\tau_q^d\Big(\frac{D}{q}\Big)e\Big(\frac{-\sum_{j}\frac{c_jc_j^{\prime}}{2\alpha_j}}{q}\Big){\sum_{a}}^* \Big(\frac{a}{q} \Big)^d e\Big(\frac{a(t-\sum_{j}\frac{c_j^{\prime 2}}{4\alpha_j})-a^{-1}\sum_{j}(\frac{c_j^2}{4\alpha_j})}{q}   \Big).$$
Next, we analyze the sum over $a$.  We define 
$$Kl(G,\vec{c},\vec{c}^{\prime},t):=e\Big(\frac{-\sum_{j}\frac{c_jc_j^{\prime}}{2\alpha_j}}{q}\Big){\sum_{a}}^* \Big(\frac{a}{q} \Big)^d e\Big(\frac{a(t-\sum_{j}\frac{c_j^{\prime 2}}{4\alpha_j})-a^{-1}\sum_{j}(\frac{c_j^2}{4\alpha_j})}{q}   \Big).$$
Note that $\Big(\frac{a}{q} \Big)^d=\Big(\frac{a}{q} \Big)$ for odd $d$. So,  $Kl(G,\vec{c},\vec{c}^{\prime},t)$ is a Sali\'e sum for odd $d$ and from standard bounds on Sali\'e sums we have
\begin{equation}\label{bound1} |Kl(G,\vec{c},\vec{c}^{\prime},t)| \leq \tau(q)\sqrt{q}.\end{equation}
Similarly, $\Big(\frac{a}{q} \Big)^d=1$  for even $d$. So, $Kl(G,\vec{c},\vec{c}^{\prime},t)$ is a Kloosterman sum for even $d$ and by Weil's bound on Kloosterman sums we obtain 
\begin{equation}\label{bound2}Kl(G,\vec{c},\vec{c}^{\prime},t) \leq \tau(q)\sqrt{q} \gcd \Big(q,(t-\sum_{j}\frac{c_j^{\prime 2}}{4\alpha_j}), \sum_{j}(\frac{c_j^2}{4\alpha_j}) \Big)^{1/2}.\end{equation}
This concludes the first part of our lemma. Finally, we analyze $S_1$.
We note that by a change of variables $S_1=S(G,\vec{c},\vec{c}^{\prime},t)$, where  $G=(mq_2)^2F,$ $t=-q_2k_1,$ and $\vec{c}^{\prime}=2q_2A\lambda.$
Recall that $F(X)=X^{T}AX$,  $\gcd(q_1,2m q_2D)=1$, and $G=(mq_2)^2F$ is diagonalizable with eigenvalues $\{\alpha_i \}$, so that
 $$\sum_{j}\frac{c_j^{\prime 2}}{4\alpha_j} \equiv  (2mq_2)^{-2} c^{\prime T} A^{-1} c^{\prime} \text{ mod }  q_1.$$
 We substitute $\vec{c}^{\prime}=2q_2A\lambda$ and obtain 
 \begin{equation} \label{subcong} 
 \sum_{j}\frac{c_j^{\prime 2}}{4\alpha_j} \equiv  \frac{\lambda^{T}A\lambda}{m^2} \equiv \frac{F(\lambda)}{m^2} \text{ mod } q_1.
 \end{equation}
 We apply formula (\ref{formulaaw}) and obtain
$|S_1|=q^{d/2}_1 |Kl(G,\vec{c},\vec{c}^{\prime},t)|.$
If $d$ is odd then by inequality (\ref{bound1}) we obtain
$|S_1|\leq  q_{1}^{\frac{d+1}{2}}\tau(q_1).$
If $d$ is even then by inequality (\ref{bound2}) we obtain
$|S_1|\leq  q_{1}^{\frac{d+1}{2}}\tau(q_1) \gcd (q_1, t-\sum_{j}\frac{c_j^{\prime 2}}{4\alpha_j})^{1/2}.$
We substitute $t=-q_2k_1$ and  by  (\ref{subcong})  obtain 
$$|S_1|\leq  q_{1}^{\frac{d+1}{2}}\tau(q_1) \gcd \big(q_1, q_2k_1+  \frac{F(\lambda)}{m^2}\big)^{1/2}.$$
By (\ref{decomk})
and (\ref{defk}), we have
 $\gcd (q_1, q_2k_1+  \frac{F(\lambda)}{m^2})= \gcd(q_1,N).$ Hence,
$$|S_1|\leq q_{1}^{\frac{d+1}{2}}\tau(q_1)\gcd(q_1,N)^{1/2}.$$
This concludes our lemma. \end{proof}

In this lemma we give an upper bound on $S_2$. This lemma is a variant of \cite[Lemma 25]{Heat} and its proof follows from the standard square root cancelation in the Gaussian sums. 
\begin{lem}\label{2} We have \begin{equation*}
|S_2|\ll_{\Delta} m^dq_2^{1+d/2}.
\end{equation*}
\end{lem}
\begin{proof}
By the Cauchy-Schwarz inequality on the $a_2$ variable, we have
\begin{equation}\begin{split}|S_2|^2&\leq m^2 \phi(q_2){\sum_{a_2}}^* \Big| \sum_{\vec{b}_2}    e\Big(\frac{2q_1a_2\lambda^TA\vec{b}_2+a_2mq_1^2F(\vec{b}_2)+\left<\vec{c},\vec{b}_2\right>-q_1a_2k_2}{mq_2} \Big)    \Big|      ^2
\\
&= m^2\phi(q_2){\sum_{a_2}}^*\sum_{\vec{b}_2 , \vec{b}_2^{\prime} }    e\Big(\frac{2q_1a_2\lambda^TA(\vec{b}_2-\vec{b}_2^{\prime})+a_2mq_1^2(F(\vec{b}_2)-F(\vec{b}_2^{\prime}))+\left<\vec{c},\vec{b}_2-\vec{b}_2^{\prime}\right>}{mq_2}     \Big),
\end{split}
\end{equation}
where $m|2\lambda^TA\vec{b}_2-k_2$ and  $m|2\lambda^TA\vec{b}^{\prime}_2-k_2$. We change the variables and write
$\vec{u}=\vec{b}_2-\vec{b}_2^{\prime}$.  Hence, 
$$|S_2|^2\leq m^2\phi(q_2){\sum_{a_2}}^*\sum_{\vec{u} , \vec{b}_2 }    e\Big(\frac{2q_1a_2\lambda^TA\vec{u}+a_2mq_1^2(2\vec{b}_2A\vec{u}+F(\vec{u}))+\left<\vec{c},\vec{u}\right>}{mq_2}     \Big),$$
where $m|2\lambda^TA\vec{b}_2-k_2$ and $m|2\lambda^TA\vec{u}$. The summation over $\vec{b}_2$ is zero unless $q_2|\Delta\gcd(\vec{u}).$
In other words, the summation is non-zero only if $$\vec{u} \in \big(q_2\mathbb{Z}/\gcd(\Delta,q_2)mq_2\mathbb{Z}\big)^d\equiv  \big(\mathbb{Z}/\gcd(\Delta,q_2) m \mathbb{Z}\big)^d.$$
Since  $\vec{b}_2 \in \big(\mathbb{Z}/mq_2\mathbb{Z}\big)^d$, then
\begin{equation*}
\begin{split}
|S_2|^2 &\leq m^2\phi(q_2){\sum_{a_2}}^*\sum_{u \in \big(\mathbb{Z}/\gcd(\Delta,q_2) m \mathbb{Z}\big)^d }\sum_{\vec{b}_2 \in \big(\mathbb{Z}/mq_2\mathbb{Z}\big)^d} 1  
\\
&= m^2 \phi(q_2)^2 \Delta^d m^{2d-2}q_2^d,
\end{split}
\end{equation*}
where $m|2\lambda^TA\vec{b}_2-k_2$ and $m|2\lambda^TA\vec{u}$. Hence,
$$|S_2| \ll_{\Delta} m^dq_2^{1+\frac{d}{2}}.$$
This concludes the lemma.
\end{proof}
\begin{rem}
We note that the exponent $1+\frac{d}{2}$ comes from the square root cancellation in the $b$ variable, and after the Cauchy-Schwarz inequality in the first line of the proof we lose all the potential cancellation in the  $a$ variable. More precisely, the sum over $a$ in the second line of the proof is over the square norm of some $d$-dimensional Gauss sums over $b$. It follows from our proof that for fixed $a$ each Gauss sum is at least $m^{d-1}q_2^{d/2}$ up to a constant which only depends on the discriminant $\Delta.$ This implies  our upper bound in Lemma~\ref{2} is sharp after the Cauchy-Schwarz inequality in the second line of the proof.  
\end{rem}

%
%

\begin{proof}[Proof of Proposition~\ref{mM}]This is a consequence of Lemma~\ref{1} and Lemma~\ref{2}. We factor $q=q_1q_2,$
where $\gcd(q_1,2m\Delta)=1$ and the prime factors of $q_2$ are a subset of the prime factors of $m\Delta $. From the formula~(\ref{FF}), we deduce that 
$|S_{m,q}(\vec{c})|=|S_1||S_2|.$
By Lemma~\ref{1}  and Lemma~\ref{2}, 
%
we have
$$\sum_{q=1}^{X}m^{-d}q^{-\frac{d+1}{2}}|S_{m,q}(\vec{c})| \ll \sum_{q=1}^{X} q_1^{\epsilon}q_2^{1/2} \gcd(q_1,N)^{1/2}
\leq X^{\epsilon}\sum_{q=1}^{X} q_2^{1/2} \gcd(q_1,N)^{1/2}.$$
 By the Cauchy-Schwarz inequality,
 $$
 \sum_{q=1}^{X} q_2^{1/2} \gcd(q_1,N)^{1/2}\leq  \big(\sum_{q=1}^{X} q_2 \big)^{1/2} \big(\sum_{q=1}^{X} \gcd(q_1,N)\big)^{1/2}.
 $$
It is easy to check that $ \sum_{q=1}^{X} \gcd(q_1,N) \leq X^{1+\epsilon}.$ For the other term, we have 
\begin{equation*}
\begin{split} \sum_{q=1}^{X} q_2 &\leq \sum_{d|(2m\Delta)^{\infty}, d<X} d \lfloor\frac{X}{d}\rfloor \leq X\sum_{d|(2m\Delta)^{\infty}, d<X} 1
 \leq X  (m\Delta X)^{\epsilon}.
\end{split}
\end{equation*}
Therefore, we conclude the proposition. \end{proof}

\section{ Construction of the smooth weight function $w$}\label{weightfunction}
In this section we construct $w\in C_{c}^{\infty}(\mathbb{R}^n)$  that we use in the delta method. 
We define   $w$ in an appropriate coordinate system, which makes the computations easier for the oscillatory integrals $I_{m,q}(\vec{c})$. 
Finally we give upper bounds on the partial derivatives of $w$, that we use in the next section.

We introduce the coordinate system first. Let 
\begin{equation}\label{vecv}
\vec{v}:=\sqrt{N}\vec{a_{\infty}}\in V_N(\mathbb{R}),
\end{equation}
 and let $e_d$ be the norm 1 vector in the direction of $\nabla F_{|\vec{v}}$ (the gradient of  $F$ at $\vec{v}$ that is $2A \vec{v}$). Let $T_{\vec{v}}(V_N(\mathbb{R}))$ be  the tangent space of $V_N(\mathbb{R})$ at $\vec{v}$, which is the orthogonal complement of $e_d$.
We restrict the quadratic from $F(\vec{x})$ to the $(d-1)$-dimensional subspace $T_{\vec{v}}(V_N(\mathbb{R}))$. By a standard theorem in Linear Algebra, we can find an orthogonal basis $B_1:=\{e_1,\dots,e_{d-1} \} $
for $T_{\vec{v}}(V_N(\mathbb{R}))$, such that 
$$
F\big(\sum_{i=1}^{d-1} u_i e_i\big)= \sum_{i=1}^{d-1} \mu_i u_i^2.
$$
Next, we construct the smooth weight function $w$ that satisfies the  required  conditions~(\ref{property1}) and~(\ref{property2}).
By these conditions, it follows that the smooth weight function $w$ is supported inside a cylinder centered at $\vec{v}$ with height $r^2\sqrt{N}$ and  radius $r\sqrt{N}$, such that its  base is parallel to $T_{\vec{v}}(V_N(\mathbb{R}))$. More precisely, let 
$$y(\vec{x}):=\frac{F(\vec{x})-N}{m^2Q^2}.$$  
Since $F(\vec{v})=\vec{v}^TA\vec{v}=N\neq0$ and $e_i^TA\vec{v}=0$ for $1\leq i\leq d-1$, then 
\begin{equation}\label{B2}B_2:=\{e_1,\dots, e_{d-1},\vec{v}\}
\text{ and }B_3:=B_1\cup\{e_d \}
\end{equation}
 are  bases for $\mathbb{R}^d$. Given $\vec{x}\in \mathbb{R}^d$, we express the vector $\vec{x}-\vec{v}$ in basis $B_3$ as $\vec{x}-\vec{v}=\sum_{i=1}^{d-1} m u_i e_i +m\alpha e_d,$ and in $B_2$ as follows:
\begin{equation}\label{corsys}\vec{x}-\vec{v}=\sum_{i=1}^{d-1} m \tilde{u_i} e_i +m\beta \vec{v}.\end{equation}
 If $\vec{v}= \sum_{i=1}^{d} v_i e_i$ then
$u_i=\tilde{u_i}+\beta v_i,$
 and $\alpha=\beta v_d.$
 Suppose $\psi_1 \in C^{\infty}_c(\mathbb{R})$ and $\psi_2 \in C^{\infty}_c(\mathbb{R}^{d-1})$. We define the smooth weight function $w(\vec{x})$ as follows:
\begin{equation}\label{weightfun}w(\vec{x})=\begin{cases}\frac{2\vec{x}^TA\vec{v}}{\sqrt{N}\left<\vec{v},e_d\right>} \psi_1(y)^2 \psi_2(\frac{\tilde{\vec{u}}}{Q})   &\text{ if  } 1+m\beta >1/2,\\
0       &\text{  otherwise,}
\end{cases}\end{equation}
where $\tilde{\vec{u}}:=(\tilde{u}_1,\dots,\tilde{u}_{d-1}).$
\begin{rem}
The factor $\frac{2\vec{x}^TA\vec{v}}{\sqrt{N}\left<\vec{v},e_d\right>}$ is in our weight function $w$  in order to simplify the oscillatory integral $I_{m,q}(\vec{c})$  in Lemma~\ref{cc}. Note that the support of the weight function $\psi_1(y)^2 \psi_2(\frac{\tilde{\vec{u}}}{Q})   $ is localized around $\vec{v}$ and $-\vec{v}$. We define $w$ such that it only has a support near $\vec{v}$.
\end{rem}
\noindent In Lemma~\ref{boundev}, we check that $w$ satisfies the required conditions~(\ref{property1}) and~(\ref{property2}). Moreover, we   give upper bounds on its partial derivatives. 
\begin{lem}\label{boundev}
Write $w=w_0\psi_1$ where $w_0(\vec{x}):=\frac{2\vec{x}^TA\vec{v}}{\sqrt{N}\left<\vec{v},e_d\right>}\psi_1(y) \psi_2(\frac{\tilde{\vec{u}}}{Q}) .$ Consider the coordinate system $(x_1,\dots,x_d)$ where $x_i=\tilde{u_i}/Q$ and $x_d=\beta N/mQ^2.$ Then  \begin{equation*} \label{supportw}w_0(\vec{x} )=0 \text{ if } \max_{1\leq i\leq d}|{x_i}| >C, \end{equation*}
 where $C$  is a constant that only depends on $\Omega$, $A$ and the support of $\psi_1$ and $\psi_2$. Moreover,  for every $1 \leq i \leq d$ and  $n \geq 0$, we have
$$|\label{boundw}\frac{\partial^n w_0}{\partial x_i^n }|<C_n,$$
 where the constant  $C_n$ only depends on $A$, compact set $\Omega$, and  $\max_{j}\big( \frac{d^j \psi_1}{d y^j}\big)$ and  $\max_{i,j} \big(\frac{\partial^j \psi_2}{\partial x_i^j }\big)$ where $0\leq j \leq n$.
\end{lem}
\begin{proof}
Since $\psi_1 \in C^{\infty}_c(\mathbb{R})$ and $\psi_2 \in C^{\infty}_c(\mathbb{R}^{d-1})$,
   there exists a constant $C^{\prime}$, such that if $\max_{1\leq i \leq d-1}|x_i|>C^{\prime}$ or  $y> C^{\prime}$ then $ w_0(\vec{x})=0$. Assuming that  $w_0(\vec{x})\neq 0$, then $1+m\beta >1/2$, $y\leq C^{\prime}$ and $|x_i|<C^{\prime}$  for $1 \leq i \leq d-1 $.  We express $y$ in terms of $(x_1,\dots,x_{d-1},x_d)$  
\begin{equation}\label{defy}
\begin{split}
y&=\frac{\big( \sum m \tilde{u}_i e_i +(1+m\beta)\vec{v} \big)^T A \big( \sum m \tilde{u}_i e_i +(1+m\beta)\vec{v} \big)-N}{m^2Q^2}
\\
&=\frac{\sum_{1}^{d-1} \tilde{u_i}^2m^2\mu_i+(1+m\beta)^2N-N}{m^2Q^2}
\\
&=\sum_{1}^{d-1} x_i^2\mu_i+ x_d(2+x_d\frac{m^2Q^2}{N}).\end{split}
\end{equation}
Therefore,
$|x_d(2+x_d\frac{m^2Q^2}{N})|\leq C^{\prime}+\sum_{1}^{d-1} C^{\prime 2}|\mu_i|.$
By \eqref{weightfun}, we have  
$$(2+x_d\frac{m^2Q^2}{N})=2+m\beta >3/2.$$ 
Then 
$ |x_d|<  C^{\prime}+\sum_{1}^{d-1} C^{\prime 2}|\mu_i|. $
We define $C:=C^{\prime}+\sum_{1}^{d-1} C^{\prime 2}|\mu_i|.$ This concludes the first part of our lemma.

Next,  we give upper bounds on the partial derivatives of $w$. We assume that $|x_i|<C$ for $1 \leq i \leq d$.    We apply  Leibniz's product rule and obtain
$$\frac{\partial^n w_0}{\partial x_i^n }= \sum_{j_1+j_2+j_3=n} {n \choose j_1, j_2, j_3}  \frac{\partial^{j_1} \frac{\vec{x}^TA\vec{v}}{\sqrt{N}\left<\vec{v},e_d\right>}}{\partial x_i^{j_1} }   \frac{\partial^{j_2}  \psi_1(y)  }{\partial x_i^{j_2} }  \frac{\partial^{j_3}  \psi_2(x_1,\dots,x_{d-1}) }{\partial x_i^{j_3} }. $$  
Hence, it suffices to show that the partial derivatives of each factor on the right hand side is $O(1)$. 
By \eqref{corsys},
\begin{equation*}
\begin{split}
\frac{\vec{x}^TA\vec{v}}{\sqrt{N}\left<\vec{v},e_d\right>}&=\frac{\vec{v}^TA\vec{v}}{\sqrt{N}\left<\vec{v},e_d\right>}+ \frac{(\vec{x}-\vec{v})^TA\vec{v}}{\sqrt{N}\left<\vec{v},e_d\right>}
\\
&=\frac{\vec{v}^TA\vec{v}}{\sqrt{N}\left<\vec{v},e_d\right>}+ \sum_{i=1}^{d-1} m Q x_i\frac{e_i^TA\vec{v}}{\sqrt{N}\left<\vec{v},e_d\right>} +(m^2Q^2x_d/N) \frac{\vec{v}^TA\vec{v}}{\sqrt{N}\left<\vec{v},e_d\right>}.
\end{split}
\end{equation*}
Since $\vec{v}^TA\vec{v}=N$ and  $e_i^TA\vec{v}=0$ for $1 \leq i \leq d-1$, 
$$\frac{\vec{x}^TA\vec{v}}{\sqrt{N}\left<\vec{v},e_d\right>}=\frac{\sqrt{N}}{\left<\vec{v},e_d\right>}+ \frac{(m^2Q^2) }{\sqrt{N}\left<\vec{v},e_d\right>} x_d.$$
We have $x_d <C$ in the support of $w$. Since $N=m^2r^{-2}Q^2$, 
$$
\frac{(m^2Q^2) }{\sqrt{N}\left<\vec{v},e_d\right>}\leq \frac{\sqrt{N}}{\left<\vec{v},e_d\right>}=\sqrt{\vec{a}_{\infty}^TA^2\vec{a}_{\infty}}=O(1).
$$ 
By the chain rule and Leibniz's product rule, we write $\frac{\partial^{j}  \psi_1(y)  }{\partial x_i^{j} }$ in terms of $ \frac{\partial^{j_1} \psi_1}{\partial y^{j_1} }$ and $\frac{\partial^{j_2} y}{\partial x_i^{j_2} }$, where $1 \leq j_1,j_2 \leq j$. Hence, it suffices to show that $\frac{\partial^j y}{\partial x_i^j }=O(1).$
By \eqref{defy},  $y$ is a quadratic form in  $(x_1,\dots,x_{d-1},x_d).$
Consequently, it suffices to show that the coefficients $\mu_i$ and $m^2Q^2/N$ are $O(1)$. Note that the $\mu_i$ are bounded by  $|A|$ (the norm of $A$) and $m^2Q^2/N<1$. Therefore, we conclude our lemma.
\end{proof}
\section{ The integral $I_{m,q}(\vec{c})$}\label{sintegral}
In this section we give upper bounds on the integral $I_{m,q}(\vec{c})$.  We use the same notations as in  Section \ref{weightfunction}. We partition the integral vectors $\vec{c}\neq 0 \in \mathbb{Z}^{d}$ into two  sets, and  denote them by ordinary and exceptional vectors. 
%
In what follows, we describe them.
   Recall  the  orthonormal basis  $B_3$ defined in  \eqref{B2}. 
We write a given integral vector $\vec{c}$ in this basis as
\begin{equation}\label{vectorc}\vec{c}=\sum_{i=1}^{d}c_ie_i.\end{equation} We remark that $\vec{c}$ does not necessarily have integral coordinates anymore. We define two types of  ordinary vectors.  The type I ordinary vectors $\vec{c}$ are the integral vectors $\vec{c}$, such that 
\begin{eqnarray}\label{typeI}
|\vec{c}|\geq (mr^{-1})^{1+\epsilon},
\end{eqnarray}
 where $|\vec{c}|=\max_{i} |c_i|$. The type II ordinary vectors $\vec{c}$ are the integral vectors $\vec{c}$, such that
\begin{equation*}
 \begin{split}
 \max_{1\leq i \leq d-1}c_i  \geq m(mr^{-1})^{\epsilon},
 \\
 |\vec{c}|< (mr^{-1})^{1+\epsilon}.
 \end{split}
\end{equation*}
We call the complement of these integral vectors the exceptional integral vectors. 

Recall that $\vec{x}=m\vec{t}+\lambda$ where $\vec{t}=(t_1,\dots,t_d)$ and let $dt=dt_1\dots dt_d$. We record some formulas for the volume form $dt$ in different coordinate systems. Recall $B_2$  from \eqref{B2} and its coordinate system $(\tilde{u_1},\dots ,\tilde{u}_{d-1},\beta)$. We have
\begin{equation}\label{change1}dt=\left<\vec{v},e_d\right>d\tilde{u_1}\dots d\tilde{u}_{d-1}d\beta.\end{equation}
%
Next, we give a formula for $dt$ in terms of $d\tilde{u}_1\dots d\tilde{u}_{d-1}dy$. By \eqref{defy}, we have  $\frac{\partial y}{\partial \beta}=\frac{2\vec{x}^TA\vec{v}}{mQ^2}$. Hence,
\begin{equation}\label{change2}
dt=\left<\vec{v},e_d\right>d\tilde{u_1}\dots d\tilde{u}_{d-1}d\beta =\frac{mQ^2}{2\vec{x}^TA\vec{v}} d\tilde{u_1}\dots d\tilde{u}_{d-1} dy.
\end{equation}
Finally, we change the variables as in lemma~\ref{boundev} to $x_i=\tilde{u_i}/Q$ for $1\leq i\leq d-1$ and $x_d=\beta N/mQ^2$
\begin{equation}\label{changex}
dt= \frac{mQ^{d+1}\left<\vec{v},e_d\right>}{2\vec{x}^TA\vec{v}}dx_1\dots dx_{d-1}dy= \frac{mQ^{d+1}\left<\vec{v},e_d\right>}{N}dx_1\dots dx_{d-1}dx_d.
\end{equation}

\subsection{ Upper bound on $I_{m,q}(\vec{c})$ for the ordinary vectors $\vec{c}$}\noindent \label{easy}

In this section we assume that  $\vec{c}$ is an ordinary vector, and prove that $I_{m,q}(\vec{c})$  rapidly decays  as $|\vec{c}|\to \infty.$   Moreover, we show that the contribution of the ordinary vectors $\vec{c}$ in the delta method is bounded. Part of our analysis in this section  is inspired by the work of Heath-Brown. In particular, the use of Fourier inversion in the proof of the following lemma is similar to  \cite[Lemma~17]{Heat}.

	\begin{lem}\label{c}
Let $\vec{c}$ be an ordinary type I vector. Then for any $k>0$, we have 
\begin{equation}\label{type1}
I_{m,q}(\vec{c}) \ll_{k} \Big(|\vec{c}|r^{-1}m  \Big)^{-k}.
\end{equation}
 On the other hand,  if $\vec{c}$ is an ordinary type II vector, we have 
\begin{equation}\label{type22}
I_{m,q}(\vec{c}) \ll_{k} \Big(r^{-1}m  \Big)^{-k}.
\end{equation}
Moreover,   
\begin{equation*}
\frac{C_{Q}}{mQ^2}\sum_{q} {\sum_{\vec{c}}}^{\prime}(mq)^{-d}S_{m,q}(\vec{c}) I_{m,q}(\vec{c}) \ll 1,
\end{equation*}
where $ {\sum}^{\prime}$ is the sum over all ordinary vectors and the implied constant in $\ll$ only depends on $\epsilon$ in the definition of the ordinary vectors, $\Omega$, and the fixed functions $\psi_1$ and $\psi_2.$
\end{lem}

%

\begin{proof}Recall from \eqref{intt}
$$I_{m,q}(\vec{c})= \int h\Big(\frac{q}{Q},\frac{G(\vec{t})}{Q^2}\Big) w(m\vec{t}+\lambda) e\big(-\frac{\left< \vec{c},\vec{t} \right>}{mq}    \big) dt.$$
By~\eqref{weightfun}, $w(m\vec{t}+\lambda)=w_0(m\vec{t}+\lambda)\psi_1(y),$
where $y=\frac{F(m\vec{t}+\lambda)-N}{m^2Q^2}=\frac{G(\vec{t})}{Q^2}$. Let $p_{\frac{q}{Q}}(\xi)$ be the Fourier transform of $h(\frac{q}{Q},y)\psi_1(y)$, we have 
$$p_{\frac{q}{Q}}(\xi):=\int_{-\infty}^{\infty}  \psi_1(y)h(\frac{q}{Q},y)e(-\xi y) dy.$$
Therefore,
$$I_{m,q}(\vec{c})=\int_{\xi} p_{\frac{q}{Q}}(\xi)\int_{t} w_0(m\vec{t}+\lambda) e\big( \xi y-\frac{\left<\vec{t}, \vec{c}\right>}{mq}   \big) dt d\xi.$$
We change the coordinate system to $(\tilde{u}_1,\dots,\tilde{u}_{d-1},\beta)$ and apply formula~(\ref{change1}) to obtain
\begin{equation}\label{fourierexp}I_{m,q}(\vec{c})=\left<\vec{v} ,e_d\right>\int_{\xi} p_{\frac{q}{Q}}(\xi)\int_{t} w_0((\tilde{\vec{u}},\beta)) e\big( \xi y-\frac{\left<\vec{t}, \vec{c}\right>}{mq}   \big) d\tilde{u}d\beta d \xi,\end{equation}
where $d\tilde{u}=d\tilde{u_1}\dots d\tilde{u}_{d-1}$.
By \eqref{defy}, we have 
$$
y=\frac{\sum_{1}^{d-1} \tilde{u_i}^2m^2\mu_i+(1+m\beta)^2N-N}{m^2Q^2}.
$$
We also write the inner product $\left<\vec{t}, \vec{c}\right>$ in terms of $\tilde{u}_1,\dots,\tilde{u}_{d-1},\beta$, and obtain  
\begin{equation}\label{innerp}
\left<\vec{t}, \vec{c}\right>=\left<\frac{\vec{v}-\lambda}{m}, c\right> +\left<\frac{\vec{x}-\vec{v}}{m}, c\right> 
=\left<\vec{t}_0, \vec{c}\right> + \sum_{i=1}^{d-1}c_i \tilde{u}_i+\beta\left<\vec{c},\vec{v} \right>,
\end{equation}
where $\vec{v}=m\vec{t}_0+\lambda$. We substitute the above expressions in equation~({\ref{fourierexp}}), and obtain
\begin{equation*}
\begin{split}
I_{m,q}(\vec{c})&=e\big(-\frac{\left<\vec{c}, \vec{t}_0\right>}{mq}    \big)\left<\vec{v},e_d\right>\int_{\xi} p_{\frac{q}{Q}}(\xi)  
\\
  &\times\int_{\beta,\tilde{u}} w_0(\tilde{u_1},\dots, \tilde{u}_{d-1}, \beta)  e\big( \frac{\xi \beta(2/m+\beta)N}{Q^2}- \frac{\beta \left<\vec{v} ,\vec{c}\right>}{mq}\Big)
\\
&\times e\big( \sum_{i=1}^{d-1}[\frac{\xi \mu_i\tilde{u}_i^2}{q^2} -\frac{c_i\tilde{u}_i}{mq}   ]\big) d\beta d\tilde{u}d\xi.
\end{split}
\end{equation*}
We are only concerned with $|I_{m,q}(\vec{c})|$, so we drop  $e\big(-\frac{\left<\vec{c}, \vec{t}_0\right>}{mq}    \big)$. We change the variables to $x_i=\tilde{u_i}/Q$ and $x_d=\beta N/mQ^2.$ By equation~(\ref{changex}), we obtain 
\begin{equation*}
\begin{split}
|I_{m,q}(\vec{c})|&=\frac{Q^{d+1}m \left<\vec{v},e_d\right>}{N} \Big|\int_{\xi} p_{\frac{q}{Q}}(\xi)  
\\
 &\int w_0(\vec{x})  e\Big( x_d(2\xi+\frac{q^2\left<\vec{v} ,\vec{c}\right>}{qN}) +\frac{\xi x_d^2m^2Q^2}{N}\Big)
\\
& e\big( \sum_{i=1}^{d-1}[\xi \mu_ix_i^2 +\frac{c_iQx_i}{mq}   ]\big) dx_1\dots dx_{d} d\xi\Big|.
\end{split}
\end{equation*}
At this point, we assume that $\vec{c}$ is an ordinary vector of type II, which means there exist $1\leq i \leq d-1$ such that  
$ c_i  \geq m(mr^{-1})^{\epsilon}.$
We write 
$$I_{m,q}(\vec{c})=I_{m,q}(\vec{c})_1+I_{m,q}(\vec{c})_2,$$ 
where $I_{m,q}(\vec{c})_1$ is the integral over  $|\xi|< (r^{-1}m)^{\epsilon/2}\frac{Q}{q}$, and $I_{m,q}(\vec{c})_2$ is the integral over $|\xi|> (r^{-1}m)^{\epsilon/2}\frac{Q}{q}$. First,  we show that 
\begin{equation}\label{type2}
| I_{m,q}(\vec{c})_2|\ll (r^{-1}m)^{-k}.
\end{equation}
Lemma \ref{Fourier} implies that  $p_{\frac{q}{Q}}(\xi)\ll(|\xi|\frac{q}{Q})^{-k}.$
Since $|\xi|>(r^{-1}m)^{\epsilon/2}\frac{Q}{q}$, we deduce the claimed upper bound on $|I_{m,q}(\vec{c})_2|$. It remains to show that 
$|I_{m,q}(\vec{c})_1|\ll (r^{-1}m)^{-k}.$
We first take the integration on the $x_i$ variable and obtain
\begin{equation*}
\begin{split}
|I_{m,q}(\vec{c})|&=
\frac{Q^{d+1}m \left<\vec{v},e_d\right>}{N} \Big|\int_{\xi} p_{\frac{q}{Q}}(\xi)  
\\
 &\int  e\Big( x_d(2\xi+\frac{q^2\left<\vec{v} ,\vec{c}\right>}{qN}) +\frac{\xi x_d^2m^2Q^2}{N}\Big)
\\
& e\big( \sum_{i=1}^{d-1}[\xi \mu_ix_i^2 +\frac{c_iqx_i}{mq}   ]\big) dx_1\dots \hat{dx}_i \dots dx_{d} d\xi
\\
&\int w_0(x_0,\dots ,x_i, \dots,x_d )e\big( \xi \mu_ix_i^2 +\frac{c_iQx_i}{mq}   \big) dx_i \Big|.\end{split}
\end{equation*}
By Lemma~\ref{boundev}, the weight function $w_0(x_0,\dots ,x_i, \dots,x_d )$ has compact support in a fixed interval $[-C,C]$, and its partial derivatives are bounded  $|\frac{\partial^n w_0}{\partial x_i^n }|<C_n$.   Since $|\xi|< (r^{-1}m)^{\epsilon/2}\frac{Q}{q}$ and $ c_i  \gg m(mr^{-1})^{\epsilon}$, we deduce the following lower bound on the derivative of the phase function  in the $x_i$ variable  
\begin{equation*}
\frac{\partial( \xi \mu_ix_i^2 +\frac{c_iQx_i}{mq})}{\partial x_i}\gg \frac{c_iQx_i}{mq}\gg \frac{(mr^{-1})^{\epsilon}Q}{q}\gg (mr^{-1})^{\epsilon}.
\end{equation*}
By integrating by parts multiple times in the $x_i$ variable, we have 
\begin{equation*}
\begin{split}
\int_{x_i} w(x_i)e\big( \xi \mu_ix_i^2 +\frac{c_iqx_i}{mq}   \big) dx_i &\ll_{N} (mr^{-1})^{-k\epsilon},
\end{split}
\end{equation*}
for every $k>0$ .
%
%
Hence, if $\vec{c}$ is a type II vector, we proved the claimed upper bound on $I_{m,q}(\vec{c})_1$, and so on $I_{m,q}(\vec{c})$. 

Next, we  show that the contribution of the type II vectors  is bounded. By Proposition~\ref{mM} and~\eqref{type22},
\begin{equation*}
\begin{split}
\frac{C_{Q}}{mQ^2}\sum_{q=1}^{Q} \sum_{\vec{c}}(mq)^{-d}S_{m,q}(\vec{c}) I_{m,q}(\vec{c}) &\ll \frac{\max_{\vec{c}}{I_{m,q}(\vec{c})}}{mQ^2} \sum_{q=1}^{Q} \sum_{\vec{c}}(mq)^{-d}|S_{m,q}(\vec{c})|
\\ & \ll \frac{ \Big(r^{-1}m  \Big)^{-k}}{m} Q^{\frac{1-d}{2}+\epsilon}m^{\epsilon}\sum_{\vec{c}} 1,
\end{split}
\end{equation*}
 where $\vec{c}$ varies over type II vectors. Since $|\vec{c}|\leq (r^{-1}m)^{1+\epsilon}$,   the number of type II vectors is bounded by $(r^{-1}m)^{d(1+\epsilon)}$.  Therefore,  by choosing  $k$ large enough, we conclude our Lemma for type II vectors. 

It remains to prove the lemma for type I vectors $\vec{c}$. We write 
$$I_{m,q}(\vec{c})=I_{m,q}(\vec{c})_1+I_{m,q}(\vec{c})_2,$$ 
where  $I_{m,q}(\vec{c})_1$ is the integral over  $|\xi|< (|\vec{c}|r^{-1}m)^{\epsilon/2}\frac{q}{Q}$ and $I_{m,q}(\vec{c})_2$ is the integral over  $|\xi|> (|\vec{c}|r^{-1}m)^{\epsilon/2}\frac{q}{Q}$ (note that the definition of $I_{m,q}(\vec{c})_1$ and $I_{m,q}(\vec{c})_2$ for the  type I ordinary vectors are slightly different from the analogous ones for the type II ordinary vectors). From the same lines as in (\ref{type2}), we have
\begin{equation*}
| I_{m,q}(\vec{c})_2|\ll (|\vec{c}|r^{-1}m)^{-k}.
\end{equation*}
It remains to prove  $I_{m,q}(\vec{c})_1 \leq (|\vec{c}|r^{-1}m)^{-k}$ for the type I vectors $\vec{c}$. If $|c_d|\neq \max_{i}|c_i|$, we proceed as before. So, we assume the harder case where $|c_d|=|c|$. We  give a lower bound for  the partial derivative of the phase function in the $x_d$ variable. We have
\begin{equation*}
\begin{split}
\frac{\partial(  x_d(2\xi+\frac{Q^2\left<\vec{v} ,\vec{c}\right>}{qN}) +\frac{\xi x_d^2m^2Q^2}{N})}{\partial x_d}&\gg \frac{Q^2\left<\vec{v} ,\vec{c}\right>}{qN}-|2\xi|-|\frac{2\xi m^2Q^2}{N}|.
\end{split}
\end{equation*}
Recall that $N=Q^2m^2r^{-2}$, and $|\xi|< (|\vec{c}|r^{-1}m)^{\epsilon/2}\frac{Q}{q}$, hence
\begin{equation*}
\frac{\partial(  x_d(2\xi+\frac{Q^2\left<\vec{v} ,\vec{c}\right>}{qN}) +\frac{\xi x_d^2m^2Q^2}{N})}{\partial x_d}\gg \frac{Q^2\left<\vec{v} ,\vec{c}\right>}{qN} \gg (|\vec{c}|r^{-1}m)^{\epsilon}.
\end{equation*}
By integrating by parts multiple times, we deduce that
$
I_{m,q}(\vec{c})_1 \ll_{k} \Big(|\vec{c}|r^{-1}m  \Big)^{-k}.
$
Therefore, we concluded the claimed upper bound (\ref{type1}) for the type I vectors. In order to prove the lemma it remains to show
\begin{equation*}
\frac{C_{Q}}{mQ^2}\sum_{q=1}^{Q} \sum_{\vec{c}}(mq)^{-d}S_{m,q}(\vec{c}) I_{m,q}(\vec{c}) \ll 1,
\end{equation*}
where the sum is over the ordinary type I vectors $\vec{c}$. We note that the following sum over $\vec{c}\in \mathbb{Z}^d$  is bounded 
\begin{equation*}
\sum_{\vec{c}\neq0}\frac{1}{|\vec{c}|^{d+1}}\ll 1,
\end{equation*}
and as a result over the type I integral vectors. Hence, by choosing $k$ large enough in (\ref{type1}), we deduce our Lemma.\end{proof}

%
%
%
%
%
%
%
%
%
%

\subsection{Upper bound on $I_{m,q}(\vec{c})$ for the exceptional vectors $\vec{c}$}\noindent \label{hard}

In this section  we assume that   $\vec{c}$ is an exceptional vector, and  give an upper bound on $I_{m,q}(\vec{c})$. 
\begin{lem}\label{cc}
We have 
\begin{equation*}
\frac{|I_{m,q}(\vec{c})|}{mQ^2} \ll \frac{Q^{d-1}}{\sqrt{N}}\min\Big[{\Big(\frac{Q|\vec{c}|}{mr^{-1}q}\Big)^{-\frac{d-1}{2}},1}\Big],
\end{equation*}
where the implied constant in $\ll$ only depends on $\Omega$ and the fixed functions $\psi_1$ and $\psi_2.$
\end{lem}
\begin{proof}
%
%
%
%
%
%
%
%
%
We apply  formula (\ref{change2}) and (\ref{innerp}) and obtain 
$$\ |I_{m,q}(\vec{c})|=mQ^2  \int \frac{\left<\vec{v},e_d  \right>}{2\vec{x}^TA\vec{v}} h(\frac{q}{Q},y) w(\vec{x}) e\big(\frac{ \sum_{i=1}^{d-1}c_i \tilde{u}_i+\beta\left<\vec{c},\vec{v} \right>}{mq}    \big) d\tilde{u}_1\dots d\tilde{u}_{d-1} dy.$$
We use the weight function $w$ defined in~(\ref{weightfun}),
and obtain 
\begin{equation}\label{Imqc} \ |I_{m,q}(\vec{c})|=\frac{mQ^2}{\sqrt{N}}  \int h(\frac{q}{Q},y) \psi_1(y)^2 \psi_2(\frac{\tilde{\vec{u}}}{Q})e\big(\frac{ \sum_{i=1}^{d-1}c_i \tilde{u}_i+\beta\left<\vec{c},\vec{v} \right>}{mq}    \big) d\tilde{u}_1\dots d\tilde{u}_{d-1} dy.
\end{equation}
By identity~\eqref{defy}, we have
$$(1+m\beta)^2=\frac{N+m^2Q^2y-\sum_{i=1}^{d-1}m^2\tilde{u_i}^2\mu_i  }{N}.$$
By changing the variables to $x_i:=\tilde{u_i}/Q$ for $1\leq i\leq d-1$ and using $N=m^2r^{-2}Q^2$, we obtain
$$m\beta= \sqrt{1+r^2y-\sum_{i=1}^{d-1}r^2x_i^2\mu_i }-1.$$
By Lemma~\ref{boundev}, if $w\neq 0$ then $x_i,y<C$ for some fixed constant $C$. Hence, by writing  the Taylor series of the square root function, we have
\begin{equation}\label{beta}
\beta=\frac{1}{2}\Big[   m^{-1}r^2y- \sum_{i=1}^{d-1}m^{-1}r^2x_i^2\mu_i \Big]+\phi(x_1,\dots ,x_{d-1},y),
\end{equation}
where
 \begin{equation}
 \label{bjbj}\frac{\partial^k \phi}{\partial x_{i_1}\dots \partial x_{i_k}} =O_k\Big(r^4m^{-1}\Big),
 \end{equation}
  for every $k\geq 0$ and any point inside the support of $w$. 
We change the variables to $(x_1,\dots,x_{d-1},y)$ in  formula~(\ref{Imqc}), and  replace $\beta$ with formula~(\ref{beta}),  
\begin{equation*}
\begin{split}
  \frac{|I_{m,q}(\vec{c})|}{mQ^2}&=  \frac{Q^{d-1}}{\sqrt{N}} \Big| \int h(\frac{q}{Q},y) \psi_1(y)^2 e\Big( \frac{\left<\vec{c},\vec{v} \right>y }{2(mr^{-1})^2q}   \Big)
 \\
  &\psi_2(x_1,\dots,x_{d-1})e\Big(\frac{ \sum_{i=1}^{d-1}Qc_i x_i-\frac{m^{-1}r^2}{2}\big[ \sum_{i=1}^{d-1} x_i^2\mu_i \big]\left<\vec{c},\vec{v} \right>+\phi \left<\vec{c},\vec{v} \right>}{mq}    \Big) dx_1\dots dx_{d-1} dy\Big|.
 \end{split}
 \end{equation*}
%
%
%
We define 
\begin{equation*}
\begin{split}
\xi(x_1,\dots, x_{d-1}, y):&=\frac{ \sum_{i=1}^{d-1}Qc_i x_i-\frac{m^{-1}r^2}{2}\big[ \sum_{i=1}^{d-1} x_i^2\mu_i \big]\left<\vec{c},\vec{v} \right>+\phi \left<\vec{c},\vec{v} \right>}{mq},  
\\
L(y):&= \int \psi_2(x_1,\dots,x_{d-1})e\Big(  \xi(x_1,\dots, x_{d-1}, y) \Big) dx_1\dots dx_{d-1}.
\end{split}
\end{equation*}
By Fubini's theorem, 
$$  \frac{|I_{m,q}(\vec{c})|}{mQ^2}=  \frac{Q^{d-1}}{\sqrt{N}} \Big| \int h(\frac{q}{Q},y) \psi_1(y)^2 e\Big( \frac{\left<\vec{c},\vec{v} \right>y }{2(mr^{-1})^2q}   \Big)L(y) dy\Big| ,$$
$$ \frac{|I_{m,q}(\vec{c})|}{mQ^2}\leq  \frac{Q^{d-1}}{\sqrt{N}}  \sup_{y\in [-C,C]} \big( \psi_1(y)^2 L(y)\big) \int |h(\frac{q}{Q},y)| dy.$$
By lemma~\ref{L12bound} for $l=k=0$, we deduce that  $\int |h(\frac{q}{Q},y)| dy<C_1,$
where $C_1$ is a constant independent of $q/Q$. Recall that $\psi_1 \in C_c^{\infty}(\mathbb{R}).$
Therefore,
\begin{equation}\label{GGG}  \frac{|I_{m,q}(\vec{c})|}{mQ^2}\ll  \frac{Q^{d-1}}{\sqrt{N}}  \sup_{y\in [-C,C]} \big( L(y)\big).\end{equation}
We give an upper bound on the oscillatory integral $L(y)$. 
We have
\begin{equation}\label{111} \frac{\partial}{\partial x_i} \xi(x_1,\dots, x_{d-1}, y)=\frac{1}{mq}\Big(Qc_i- \frac{\mu_i\left<\vec{c},\vec{v} \right>x_i }{mr^{-2}}+\left<\vec{c},\vec{v} \right>\frac{\partial  \phi}{\partial x_{i} } \Big).\end{equation}
  Without loss  of generality, we assume that $r<\epsilon_0$, where $\epsilon_0$ depends only on the  compact set $\Omega$. Let $c_l=\max_{1\leq i\leq d-1} c_i$ and assume that 
   \begin{equation}\label{222}\Big|\frac{\mu_l\left<\vec{c},\vec{v} \right>|C| }{mr^{-2}} \Big|\leq \frac{Qc_{l}}{4} ,\end{equation}
   where the constant  $C$ is defined in Lemma~\ref{boundev}. By \eqref{bjbj},  $\frac{\partial  \phi}{\partial x_{l} }=O\Big(r^4m^{-1}\Big)$ in the support of $w$. Hence, by choosing $\epsilon_0$ small enough 
 \begin{equation}\label{333}\Big| \frac{\partial}{\partial x_{l} } \phi\left<\vec{c},\vec{v} \right> \Big| \leq \frac{Qc_{l}}{4}.\end{equation}
 By applying the inequalities~\eqref{222} and~\eqref{333} on equation~(\ref{111}), we obtain  
  $$ \frac{\partial}{\partial x_l} \xi \geq \frac{Qc_{l}}{2mq} .$$   
  By inequalities~\eqref{bjbj} and~\eqref{222}  for every $k\geq 0$, we have 
  $$
|\frac{\partial^k \xi}{\partial x_{l}^k}|\ll_k  \frac{Qc_l}{mq}.
  $$
  By integration by parts  multiple times in the $x_l$ variable, we deduce that 
  $$  L(y) \ll_{A}  \sup_{0 \leq j \leq A+1} |\frac{\partial^j \xi}{\partial x_{l}^j}| \Big(\frac{Q|{c_l}|}{mq}\Big)^{-(A+1)}  \ll_A\min\Big[{\Big(\frac{Q|{c_l}|}{mq}\Big)^{-A},1}\Big].$$
 Finally, by~(\ref{GGG}) and~\eqref{222}, we deduce that 
 $$\frac{|I_{m,q}(\vec{c})|}{mQ^2}\ll  \frac{Q^{d-1}}{\sqrt{N}}  \min\Big[{\Big(\frac{Q|\vec{c}|}{mr^{-1}q}\Big)^{-A},1}\Big].$$
 This concludes the lemma by assuming~(\ref{222}). It remains to prove our lemma when 
 $$\frac{Qc_{l}}{4}<\Big|\frac{\mu_l\left<\vec{c},\vec{v} \right>|C| }{mr^{-2}} \Big|,$$
for every $ 1 \leq l \leq d-1$. Since   $\vec{a_{\infty}}=\vec{v}/\sqrt{N}\in \Omega$ and $\sqrt{N}=Qmr^{-1}$ then 
$c_{l}\ll  \Big|r \left< \vec{c},\vec{a}_{\infty}  \right> \Big|.$
Let $\vec{a}_{\infty}= \sum_{i=1}^{d} a_i e_i$ then by choosing $\epsilon_0 $ small enough, we deduce that 
$$1/2 \left< \vec{c},\vec{a}_{\infty}\right> \leq  c_d a_{d} < 2\left< \vec{c},\vec{a}_{\infty}\right>, \text{ and } |c_d|= \max_{i}c_i.$$
Hence, 
\begin{equation}\label{haha}|\vec{c}| \ll  | \left< \vec{c},\vec{a}_{\infty}\right>|.\end{equation}
Therefore, 
$$ \frac{\partial^2}{\partial x_i \partial x_j} \xi(x_1,\dots, x_{d-1}, y)=\frac{1}{mq}\big(- \delta(i,j)\frac{\mu_i\left<\vec{c},\vec{v} \right>}{mr^{-2}}+\frac{\partial^2}{\partial x_{i} \partial x_{j} } \phi\left<\vec{c},\vec{v} \right> \big),$$
where $\delta(i,j)=1$ if $i=j$ and $\delta(i,j)=0$ otherwise. By applying~(\ref{bjbj}), we obtain 
$$  \frac{\partial^2}{\partial x_i \partial x_j} \xi(x_1,\dots, x_{d-1}, y)=- \delta(i,j)\frac{\mu_i\left<\vec{c},\vec{v} \right>}{m^2qr^{-2}}+O\Big(\frac{\left<\vec{c},\vec{v} \right>}{m^2qr^{-4}}   \Big).$$ 
We substitute $\sqrt{N}\vec{a}_{\infty}=\vec{v}$ and $mr^{-1}Q=\sqrt{N}$, and get
$$  \frac{\partial^2}{\partial x_i \partial x_j} \xi(x_1,\dots, x_{d-1}, y)=- \delta(i,j)\frac{\mu_i Q\left<\vec{c},\vec{a}_{\infty} \right>}{mr^{-1}q}+O\Big(\frac{Q\left<\vec{c},\vec{a}_{\infty} \right>}{qmr^{-3}}   \Big).$$ 
 Finally by inequality~(\ref{haha}) and the stationary phase theorem on the oscillatory integral $L(y)$, we obtain
\begin{equation*}\label{stationary}
L(y)\ll \frac{Q^{d-1}}{\sqrt{N}}\min\Big[{\Big(\frac{Q|\vec{c}|}{mr^{-1}q}\Big)^{-\frac{d-1}{2}},1}\Big].
\end{equation*}
This concludes the lemma.\end{proof}

 \subsection{Bounding the contribution of the exceptional vectors }
 In this section we bound the contribution of  the exceptional vectors to the right hand side of  \eqref{newequu}.
 \begin{prop}\label{cor}
We have
\begin{equation*}
\sum_{q=1}^{Q}\sum_{\vec{c}\neq 0} (mq)^{-d}S_{m,q}(\vec{c}) I_{m,q}(\vec{c}) 
\ll \begin{cases} \frac{m Q^{5}}{\sqrt{N}} (r^{-1}m)^{\epsilon}  \frac{r^{-1}m}{Q^{1/2}}  &\text{if } d=4,
\\
\\
\frac{m Q^{d+1}}{\sqrt{N}} (r^{-1}m)^{\epsilon}  \big(\frac{r^{-1}m}{Q}\big)^{\frac{d-3}{2}}   &\text{if } d\geq 5
\end{cases}
\end{equation*}
for any $\epsilon>0$, where the sum is over   the exceptional vectors $\vec{c}$, and the implied constant in $\ll$  only depends on $\epsilon$, $\Omega$ and the fixed functions $\psi_1$ and $\psi_2.$ 
\end{prop}
We give the proof of the above proposition at the end of this section. We write 
$$\sum_{q=1}^{Q}\sum_{\vec{c}\neq 0} (mq)^{-d}S_{m,q}(\vec{c}) I_{m,q}(\vec{c}) =H_1+H_2,$$
where 
$$H_1:=\sum_{\vec{c}\neq 0}\sum_{q=1 }^{\frac{Q |\vec{c}|}{r^{-1}m}} (mq)^{-d}S_{m,q}(\vec{c}) I_{m,q}(\vec{c}),$$
and 
$$H_2:=\sum_{\vec{c}\neq 0}\sum_{q=\frac{Q |\vec{c}|}{r^{-1}m} }^Q(mq)^{-d}S_{m,q}(\vec{c}) I_{m,q}(\vec{c}).$$
By Lemma~\ref{conglem}, $S_{m,q}(\vec{c})=0$, unless 
$\vec{c}\equiv \alpha A\lambda$   for  some scalar $\alpha$.  Without loss of generality we assume that $\vec{c}\equiv \alpha A\lambda   \text  {  mod }  m$.
 First, we bound  the number of exceptional vectors satisfying this congruence condition.%
 
  
 \begin{lem}\label{countt} The number of exceptional integral vectors $\vec{c}$ such that 
$\vec{c}\equiv \alpha A\lambda   \text  {  mod }  m,$
 for some $\alpha \in \mathbb{Z}/m\mathbb{Z}$  and
 $ |\vec{c}| \leq T,$
is bounded by $T(r^{-1}m)^{\epsilon}$.

 \end{lem}
\begin{proof}Fix $\alpha$ mod $m$ with one of the $m$ choices mod $m$. The proof is based on the covering of $\mathbb{R}^d$ by boxes of volume $m^{d}$ around the integral points $c$ where  $\vec{c}\equiv \alpha A\lambda$    mod  $m.$ 
By the definition of the exceptional vectors, 
$
|c_i|< m(mr^{-1})^{\epsilon}.
$
Therefore, all the exceptional vectors with $|\vec{c}|<T$ lie inside a box of volume  $Tm^{d-1}(mr^{-1})^{(d-1)\epsilon}.$
 By a covering argument,  the number of exceptional vectors is less than
 $$m\frac{Tm^{d-1}(mr^{-1})^{(d-1)\epsilon}}{m^d}=T(mr^{-1})^{(d-1)\epsilon}.$$
 By choosing $\epsilon$ small enough in the definition of the exceptional vectors, we conclude the lemma.
  \end{proof}


%
%
%
%
%

	\begin{lem}\label{XX}
We have 
\begin{equation}\label{upup}
H_1 \ll \frac{m Q^{d+1}}{\sqrt{N}}\big( \frac{Q}{r^{-1}m}\big)^{-\frac{d-3}{2}}(r^{-1}m)^{\epsilon} \big[(r^{-1}m)^{-\frac{d-5}{2}}+1 \big]
\end{equation}
for any $\epsilon>0$, where the implied constant in $\ll$  only depends on $\epsilon$, $\Omega$ and the fixed functions $\psi_1$ and $\psi_2.$ \end{lem}
%

\begin{proof}
 By Lamma~\ref{cc},
$$
H_1 \ll \frac{m Q^{d+1}}{\sqrt{N}}  \sum_{\vec{c}}\sum_{q }(mq)^{-d}|S_{m,q}(\vec{c})|\Big(\frac{Q |\vec{c}|}{r^{-1}m q}\Big)^{-\frac{d-1}{2}}.
$$
By Proposition~\ref{mM}, we have
\begin{equation}\label{sumq1}\sum_{q }m^{-d}q^{-\frac{d+1}{2}} |S_{m,q}(\vec{c})|   \ll m^{\epsilon}(\frac{Q}{r^{-1}m}|\vec{c}|)^{1+\epsilon}.\end{equation} 
By Lemma~\ref{countt}, we deduce that
\begin{equation*}
\sum_{|\vec{c}|<(r^{-1}m)^{1+\epsilon}} |\vec{c}|^{-\frac{d-3}{2}} \ll (r^{-1}m)^{\epsilon} \big[(r^{-1}m)^{-\frac{d-5}{2}}+1 \big].
\end{equation*}
Therefore,
$$
H_1 \ll\frac{m Q^{d+1}}{\sqrt{N}}\big( \frac{Q}{r^{-1}m}\big)^{-\frac{d-3}{2}}(r^{-1}m)^{\epsilon} \big[(r^{-1}m)^{-\frac{d-5}{2}}+1 \big].
$$
By choosing  $\epsilon$ small enough in the definition of the exceptional vectors $\vec{c}$, we concludes the lemma. \end{proof}
\begin{rem}
Note that if $d\geq 5$, then $(r^{-1}m)^{-\frac{d-5}{2}}=O(1)$ and the upper bound~\eqref{upup} is sharp. However, if $d=4$ by applying inequality~\eqref{sumq1} in the proof, we lose all the potential cancellation in the $q$ variable. In fact, if we assume the square root cancellation in the $q$ variable, then 
\begin{equation}\label{imo}H_1 \ll \frac{m Q^{d+1}}{\sqrt{N}}\big( \frac{Q}{r^{-1}m}\big)^{-\frac{d-3}{2}}(r^{-1}m)^{\epsilon}\end{equation}
for any $d\geq 4.$ 
\end{rem}

%
%
%
%

\begin{lem}\label{YY}
We have
\begin{equation}\label{upupup}
H_2\ll\frac{m Q^{d+1}}{\sqrt{N}} \big( \frac{Q}{r^{-1}m}\big)^{-\frac{d-3}{2}}  (r^{-1}m)^{\epsilon}  \big[(r^{-1}m)^{-\frac{d-5}{2}}+1 \big]
\end{equation}
 for any $\epsilon>0$, where the implied constant in $\ll$  only depends on $\epsilon$, $\Omega$ and the fixed functions $\psi_1$ and $\psi_2.$ \end{lem}

\begin{proof}  
By using  the trivial  bound $|I_{m,q}(\vec{c})|\ll  \frac{m Q^{d+1}}{\sqrt{N}} $ and Proposition~\ref{mM}, we have
\begin{equation}\label{sssqqq}
\begin{split}
H_2 &\ll \frac{m Q^{d+1}}{\sqrt{N}} \sum_{q}(mq)^{-d}|S_{m,q}(\vec{c})| 
\\
& \ll \frac{m Q^{d+1}}{\sqrt{N}} m^{\epsilon}\big( \frac{Q}{r^{-1}m} |\vec{c}|\big)^{-\frac{d-3}{2}+\epsilon}.
\end{split}
\end{equation}
By   Lemma~\ref{countt} ,  we obtain
$$\sum_{|\vec{c}|<(r^{-1}m)^{1+\epsilon}} |\vec{c}|^{-\frac{d-3}{2}} \ll (r^{-1}m)^{\epsilon} \big[(r^{-1}m)^{-\frac{d-5}{2}}+1 \big].$$
Hence, 
\begin{equation}
H_2\ll \frac{m Q^{d+1}}{\sqrt{N}} \big( \frac{Q}{r^{-1}m}\big)^{-\frac{d-3}{2}}  (r^{-1}m)^{\epsilon}  \big[(r^{-1}m)^{-\frac{d-5}{2}}+1 \big].
\end{equation}
This concludes the lemma.\end{proof}
\begin{rem}\label{evidence}
Similarly, if $d\geq 5$ then the upper bound~\eqref{upupup} is sharp, and if  $d=4$ by applying inequality~\eqref{sssqqq} in the proof we lose a square root cancelation  in the $q$ variable. In fact, if we assume the square root cancellation in the $q$ variable, then 
$$H_2\ll \frac{m Q^{d+1}}{\sqrt{N}}\big( \frac{Q}{r^{-1}m}\big)^{-\frac{d-3}{2}}(r^{-1}m)^{\epsilon}$$
for any $d\geq 4$. The above inequality and  inequality~\eqref{imo}  implies Conjecture~\ref{mainconjecture}. 
\end{rem}

\begin{proof}[Proof of Proposition \ref{cor}]
By Lemma~\ref{XX}, it follows that 
\begin{equation*}
H_1
\ll \begin{cases} \frac{m Q^{d+1}}{\sqrt{N}} (r^{-1}m)^{\epsilon}  \frac{r^{-1}m}{Q^{1/2}}  &\text{if } d=4,
\\
\\
\frac{m Q^{d+1}}{\sqrt{N}} (r^{-1}m)^{\epsilon}  \big(\frac{r^{-1}m}{Q}\big)^{\frac{d-3}{2}}   &\text{if } d\geq 5.
\end{cases}
\end{equation*}
Similarly, by Lemma~\ref{YY}, we have 
\begin{equation*}
H_2
\ll \begin{cases} \frac{m Q^{d+1}}{\sqrt{N}} (r^{-1}m)^{\epsilon}  \frac{r^{-1}m}{Q^{1/2}}  &\text{if } d=4,
\\
\\
\frac{m Q^{d+1}}{\sqrt{N}} (r^{-1}m)^{\epsilon}  \big(\frac{r^{-1}m}{Q}\big)^{\frac{d-3}{2}}   &\text{if } d\geq 5.
\end{cases}
\end{equation*}
The proposition follows from the above inequalities. \end{proof}

%
%
%
%

%
%
%
%
%
%

\section {The main theorem }\label{final}

Recall  the notations used while formulating Theorem~\ref{strong}, and the definition of $N(w,\lambda )$ in \eqref{newequu}. In this section we give an asymptotic formula for $N(w,\lambda )$, which implies Theorem~\ref{strong}. The main term of this formula comes from $\vec{c}=0$. 

First, we estimate the integral $I_{m,q}(0)$. Let $(x_1,\dots,x_d)$ be the coordinate system defined in Lemma~\ref{boundev}. Define 
$$
\sigma_{\infty}(F,w)=  \lim_{\epsilon \to 0}\frac{\int_{-\epsilon/2<y(x)<\epsilon/2} w(\vec{x}) \left<\vec{a}_{\infty},e_d\right>dx_1\dots dx_{d}}{\epsilon}.
$$
By  Lemma~\ref{boundev}, it follows that there exists constants $c_1$ and $c_2$, that only depend  on $\Omega$ and the fixed smooth functions $\psi_1$ and $\psi_2$, such that 
\begin{equation}\label{sigma}c_1\leq \sigma_{\infty}(F,w)\leq c_2.\end{equation}

\begin{lem}\label{singg}
Suppose that $1\leq q \leq Q^{1-\epsilon}$. Then
\begin{equation*}
I_{m,q}(0)=\sigma_{\infty}(F,w)\frac{m Q^{d+1}}{\sqrt{N}}+O_N(Q^{-k})
\end{equation*}
for any $k>0$.
\end{lem}
\begin{proof}
Recall the formula (\ref{intt}), and plug in $\vec{c}=0$ in that  to obtain
$$I_{m,q}(0)= \int h\Big(\frac{q}{Q},\frac{F(m\vec{t}+\lambda)-N}{m^2Q^2}\Big) w(m\vec{t}+\lambda) dt_1\dots dt_d.
$$
Recall the definition of $w$ from \eqref{weightfun}.
 Change the variables from $(t_1, \dots t_d)$ to $(x_1,\dots,x_{d-1},y)$ that were introduced in Lemma~\ref{boundev}, and apply \eqref{changex} to get
\begin{equation*}I_{m,q}(0)= \int \frac{mQ^{d+1}}{N} h(\frac{q}{Q},y) \psi_1(y)^2 \psi_2(x_1,\dots,x_{d-1})   dx_1\dots dx_{d-1} dy.
\end{equation*}
By taking the integration over $y$, and using Lemma \ref{del}, we have  
\begin{equation*}
\begin{split}
I_{m,q}(0)= \frac{mQ^{d+1}}{\sqrt{N}}  \int \psi_2(x_1,\dots,x_{d-1}) \psi_1(0)^2 dx_1\dots dx_{d-1}+O_N(Q^{-k}).
\end{split}
\end{equation*}
By \eqref{weightfun},  $$ \psi_2(x_1,\dots,x_{d-1}) \psi_1(0)^2=\frac{\sqrt{N}\left<\vec{v},e_d\right>} {2\vec{x}^TA\vec{v}} w(\vec{x})$$ for $y=0$. By \eqref{defy},  $\big(\frac{\partial y}{\partial x_d}\big)^{-1}=\frac{N}{2\vec{x}^TA\vec{v}}.$
We substitute these values and obtain 
\begin{equation*}
\begin{split}
I_{m,q}(0)&=
\frac{mQ^{d+1}}{\sqrt{N}} \int_{y(x)=0} \big(\frac{\partial y}{\partial x_d}\big)^{-1} w(\vec{x}) \left<\vec{a}_{\infty},e_d\right>dx_1\dots dx_{d-1}+O_N(Q^{-k})
\\
&= \frac{mQ^{d+1}}{\sqrt{N}} \lim_{\epsilon \to 0}\frac{\int_{-\epsilon/2<y(x)<\epsilon/2} w(\vec{x})\left<\vec{a}_{\infty},e_d\right> dx_1\dots dx_{d}}{\epsilon}+O_N(Q^{-k})
\\
&=\sigma_{\infty}(F,w)\frac{m Q^{d+1}}{\sqrt{N}}+O_N(Q^{-k}).
\end{split}
\end{equation*}
This concludes the proof of our lemma. 
\end{proof}
%
%
In the following lemma, we show that the contribution of $Q^{1-\epsilon}\leq q \leq Q$ is negligible in \eqref{newequu}.
\begin{lem}\label{c=0}
We have
\begin{equation*}
\sum_{q=Q^{1-\epsilon}}^{Q}(mq)^{-d}S_{m,q}(0)I_{m,q}(0)\ll_{\epsilon} \frac{m Q^{d+1}}{\sqrt{N}} Q^{\frac{3-d}{2}+\epsilon}.
\end{equation*}
\end{lem}
\begin{proof} We have  $I_{m,q}(0) \ll \frac{m Q^{d+1}}{\sqrt{N}}.$
By Proposition~\ref{mM}, we have
\begin{equation*}
\sum_{q=Q^{1-\epsilon}}^{Q}(mq)^{-d}|S_{m,q}(0)|\ll Q^{\frac{3-d}{2}+\epsilon}.
\end{equation*}
Hence, we conclude the Lemma. \end{proof}

We invoke~\cite[Page 50; Lemma 31]{Heat}.
\begin{lem}[Heath-Brown]\label{singser}
We have
\begin{equation*}
\sum_{1 \leq q \leq T} (qm)^{-d}S_{m,q}(0)=\mathfrak{S}_{B_{\vec{a},\vec{r}}}(N)+O(T^{(3-d)/2+\epsilon}Q^{\epsilon})
\end{equation*}
for any $\epsilon>0$. As a result 
\begin{equation*}
\sum_{1 \leq q \leq Q^{1-\epsilon}} (qm)^{-d}S_{m,q}(0)=\mathfrak{S}_{B_{\vec{a},\vec{r}}}(N)+O(Q^{(3-d)/2+\epsilon})
\end{equation*}
for any $\epsilon>0$, where the implied constant in $O_{\epsilon}$ only depends on $\epsilon$ and the quadratic from $F$.
\end{lem}
%
%

\begin{proof}[Proof of Theorem~\ref{strong}]
Let $w$ be the smooth weight function defined in~\eqref{weightfun}.  By Lemma~\ref{boundev},  $\text{supp}(w) \subset B_{\infty}(\vec{a}_{\infty},r).$ Recall~\eqref{newequu} 
\begin{equation*}
N(w,\lambda )=\frac{C_Q}{mQ^2}\sum_{q=1}^{Q}\sum_{\vec{c}}(mq)^{-d}S_{m,q}(\vec{c})I_{m,q}(\vec{c}).
\end{equation*}
From Proposition~\ref{cor} and Lemma \ref{c=0}, we can drop the nonzero integral vectors $\vec{c}$ and restrict the sum over $1\leq q \leq Q^{1-\epsilon}$, hence for $d\geq 5$
\begin{equation*}
N(w,\lambda)=\frac{C_Q}{mQ^2}\sum_{q=1}^{Q^{1-\epsilon}}(mq)^{-d}S_{m,q}(0)I_{m,q}(0)+O_{\epsilon}\big(|B_{\vec{a},\vec{r}}|^{d-1}N^{\frac{d-2}{2}}(|B_{\vec{a},\vec{r}}|^2\sqrt{N} )^{-\frac{d-3}{2}}N^{\epsilon}\big),
\end{equation*}
and for $d=4$, we have 
\begin{equation*}
N(w,\lambda)=\frac{C_Q}{mQ^2}\sum_{q=1}^{Q^{1-\epsilon}}(mq)^{-d}S_{m,q}(0)I_{m,q}(0)+O_{\epsilon}\big(|B_{\vec{a},\vec{r}}|^3N(|B_{\vec{a},\vec{r}}|^3\sqrt{N} )^{-\frac{1}{2}}N^{\epsilon}\big).
\end{equation*}
By Lemma~\ref{singg}, we have 
$
I_{m,q}(0)=\sigma_{\infty}(F,w)\frac{m Q^{d+1}}{\sqrt{N}}+O_{k}(Q^{-k})
$
for any $k>0.$
Therefore, 
$$\frac{C_Q}{mQ^2}\sum_{q=1}^{Q^{1-\epsilon}}(mq)^{-d}S_{m,q}(0)I_{m,q}(0)=  \sigma_{\infty}(F,w)\frac{ Q^{d-1}}{\sqrt{N}}\sum_{q=1}^{Q^{1-\epsilon}}(mq)^{-d}S_{m,q}(0)+O_{k}(Q^{-k}).$$
Finally, from  Lemma \ref{singser}, if $d\geq 5$ it follows that 
\begin{equation*}
N(w,\lambda)= \sigma_{\infty}(F,w)\mathfrak{S}_{B_{\vec{a},\vec{r}}}(N) |B_{\vec{a},\vec{r}}|^{d-1}N^{\frac{d-2}{2}}\Big(1+O_{\epsilon}\big((|B_{\vec{a},\vec{r}}|^2\sqrt{N} )^{-\frac{d-3}{2}}N^{\epsilon}\big) \Big),
\end{equation*}
and for $d=4$, we have  
\begin{equation*}
N(w,\lambda)=\sigma_{\infty}(F,w)\mathfrak{S}_{B_{\vec{a},\vec{r}}}(N) |B_{\vec{a},\vec{r}}|^{3}N\Big(1+O_{\epsilon}\big((|B_{\vec{a},\vec{r}}|^3\sqrt{N} )^{-\frac{1}{2}}N^{\epsilon}\big) \Big).
\end{equation*}
Since $w$ is bounded and supported inside $B_{\infty}(\vec{a}_{\infty},r)$, then 
$|V_N(\mathbb{Z})\cap B| \gg  N(w,\lambda).$
This concludes Theorem~\ref{strong}.
 \end{proof}
\begin{rem}
Since we derived a smooth counting formula for the number of integral points with a power saving error term. By a standard method in analysis (majorants and minorants for
the indicator function of a ball) we can find two smooth weight function $w_1(x)$ and $w_2(x)$ such that they approximate from below and above the indicator function $\chi(x)$ of the open ball $B_{\infty}(\vec{a}_{\infty},r)$, i.e. $
w_1(x) \leq \chi(x) \leq w_2(x).
$
In this way it would be possible to give a counting formula for Theorem~\ref{strong} with a power saving error term instead of the stated  lower bound.
\end{rem}

%
%
%
%

%
%

\bibliographystyle{alpha}
\bibliography{revised}

\begin{thebibliography}{{Sar}17}

\bibitem[AC37]{Chowla}
F.~C. Auluck and S.~Chowla.
\newblock The representation of a large number as a sum of ‘almost equal’
  squares,.
\newblock {\em Proc. Indian Acad. Sci.}, 6:81–82., 1937.

\bibitem[BR12]{Bourgain}
Jean Bourgain and Ze{\'e}v Rudnick.
\newblock Restriction of toral eigenfunctions to hypersurfaces and nodal sets.
\newblock {\em Geom. Funct. Anal.}, 22(4):878--937, 2012.

\bibitem[Chi95]{Chiu}
Patrick Chiu.
\newblock Covering with {H}ecke points.
\newblock {\em J. Number Theory}, 53(1):25--44, 1995.

\bibitem[Dae10]{Lower}
Dirk Daemen.
\newblock Localized solutions in {W}aring's problem: the lower bound.
\newblock {\em Acta Arith.}, 142(2):129--143, 2010.

\bibitem[DFI93]{delta}
W.~Duke, J.~Friedlander, and H.~Iwaniec.
\newblock Bounds for automorphic {$L$}-functions.
\newblock {\em Invent. Math.}, 112(1):1--8, 1993.

\bibitem[Eic54]{Eichler}
Martin Eichler.
\newblock Quatern\"are quadratische {F}ormen und die {R}iemannsche {V}ermutung
  f\"ur die {K}ongruenzzetafunktion.
\newblock {\em Arch. Math.}, 5:355--366, 1954.

\bibitem[GGN13]{GGN}
Anish Ghosh, Alexander Gorodnik, and Amos Nevo.
\newblock Diophantine approximation and automorphic spectrum.
\newblock {\em Int. Math. Res. Not. IMRN}, (21):5002--5058, 2013.

\bibitem[GGN15]{GGN1}
Anish Ghosh, Alexander Gorodnik, and Amos Nevo.
\newblock Diophantine approximation exponents on homogeneous varieties.
\newblock In {\em Recent trends in ergodic theory and dynamical systems},
  volume 631 of {\em Contemp. Math.}, pages 181--200. Amer. Math. Soc.,
  Providence, RI, 2015.

\bibitem[GGN16]{GGN2}
Anish Ghosh, Alexander Gorodnik, and Amos Nevo.
\newblock Best possible rates of distribution of dense lattice orbits in
  homogeneous spaces.
\newblock {\em J. reine angew. Math.}, 2016.

\bibitem[Har90]{Harman}
Glyn Harman.
\newblock Approximation of real matrices by integral matrices.
\newblock {\em J. Number Theory}, 34(1):63--81, 1990.

\bibitem[HB96]{Heat}
D.~R. Heath-Brown.
\newblock A new form of the circle method, and its application to quadratic
  forms.
\newblock {\em J. Reine Angew. Math.}, 481:149--206, 1996.

\bibitem[Hoo78]{Hoole}
Christopher Hooley.
\newblock On the greatest prime factor of a cubic polynomial.
\newblock {\em J. Reine Angew. Math.}, 303/304:21--50, 1978.

\bibitem[Kim03]{Kim}
Henry~H. Kim.
\newblock Functoriality for the exterior square of {${\rm GL}_4$} and the
  symmetric fourth of {${\rm GL}_2$}.
\newblock {\em J. Amer. Math. Soc.}, 16(1):139--183, 2003.
\newblock With appendix 1 by Dinakar Ramakrishnan and appendix 2 by Kim and
  Peter Sarnak.

\bibitem[Klo27]{Kloos}
H.~D. Kloosterman.
\newblock On the representation of numbers in the form {$ax^2+by^2+cz^2+dt^2$}.
\newblock {\em Acta Math.}, 49(3-4):407--464, 1927.

\bibitem[LPS88]{Rama}
A.~Lubotzky, R.~Phillips, and P.~Sarnak.
\newblock Ramanujan graphs.
\newblock {\em Combinatorica}, 8(3):261--277, 1988.

\bibitem[Mal62]{Malyshev}
A.~V. Maly{\v{s}}ev.
\newblock On the representation of integers by positive quadratic forms.
\newblock {\em Trudy Mat. Inst. Steklov}, 65:212, 1962.

\bibitem[Mar88]{Margulis}
G.~A. Margulis.
\newblock Explicit group-theoretic constructions of combinatorial schemes and
  their applications in the construction of expanders and concentrators.
\newblock {\em Problemy Peredachi Informatsii}, 24(1):51--60, 1988.

\bibitem[RS17]{Rivin}
Igor Rivin and Naser~T. Sardari.
\newblock Quantum chaos on random cayley graphs of.
\newblock {\em Experimental Mathematics}, 0(0):1--14, 2017.

\bibitem[Sar90]{10}
Peter Sarnak.
\newblock {\em Some applications of modular forms}, volume~99 of {\em Cambridge
  Tracts in Mathematics}.
\newblock Cambridge University Press, Cambridge, 1990.

\bibitem[Sar15a]{Sarnak2}
Peter Sarnak.
\newblock {\em {Letter to Scott Aaronson and Andy Pollington on the
  Solovay-Kitaev Theorem}}, February 2015.

\bibitem[Sar15b]{Sarnak3}
Peter Sarnak.
\newblock {\em {Letter to Stephen D. Miller and Naser Talebizadeh Sardari on
  optimal strong approximation by integral points on quadrics; the case:
  $SL_2(\mathbb{Z}) \to SL_2(\frac{\mathbb{Z}}{q\mathbb{Z}})$ }}, August 2015.

\bibitem[{Sar}17]{Sardaric}
N.~T {Sardari}.
\newblock {Complexity of strong approximation on the sphere}.
\newblock {\em ArXiv e-prints}, March 2017.

\bibitem[Sar18]{Naser}
Naser~T. Sardari.
\newblock Diameter of ramanujan graphs and random cayley graphs.
\newblock {\em Combinatorica}, Aug 2018.

\bibitem[Sie67]{Siegel}
C.~L. Siegel.
\newblock {\em Lectures on quadratic forms}.
\newblock Notes by K. G. Ramanathan. Tata Institute of Fundamental Research
  Lectures on Mathematics, No. 7. Tata Institute of Fundamental Research,
  Bombay, 1967.

\bibitem[Tij86]{Tj}
R.~Tijdeman.
\newblock Approximation of real matrices by integral matrices.
\newblock {\em J. Number Theory}, 24(1):65--69, 1986.

\bibitem[Wri33]{Wright}
E.~M. Wright.
\newblock The representation of a number as a sum of five or more squares (ii).
\newblock {\em Quart. J. Math. Oxford}, 4:37--51 ; 228--232., 1933.

\bibitem[Wri37]{Wright2}
E.~M. Wright.
\newblock The representation of a number as a sum of four almost equal squares.
\newblock {\em Quart. J. Math. Oxford}, 8:278--279., 1937.

\end{thebibliography}

\end{document}